\input ./style/arxiv-general.cfg
\documentclass[aop,MSNbibl,nameyear,dvips]{arximspdf}
\makeatletter
   \@ifpackageloaded{graphicx}{}{\usepackage{graphicx}}
\makeatother
\usepackage{mathrsfs}
\usepackage{mathbh}

\doi{10.1214/15-AOP1031}
\volume{44}
\issue{4}
\pubyear{2016}
\firstpage{2643}
\lastpage{2693}
\docsubty{FLA}

\makeatletter

\newcommand{\rrvert}{\vert}

\newcommand{\rrVert}{\Vert}
\newcommand{\llvert}{\vert}
\newcommand{\llVert}{\Vert}
\renewcommand{\mid}{|}
\def\sfrac#1#2{#1/#2}
\def\vfrac#1#2{(#1)/#2}

\newcommand{\mathds}{\mathbb}

\newtheorem{theorem}{Theorem}[section]
\newproclaim{definition}[theorem]{Definition}
\newproclaim{assumption}[theorem]{Assumption}
\newproclaim{example}[theorem]{Example}
\newproclaim{remark}[theorem]{Remark}

\newcommand{\Dom}{\operatorname{Dom}}
\newcommand{\tr}{\operatorname{Tr}}
\newcommand{\supp}{\operatorname{supp}}
\newcommand{\bE}{\mathds{E}}
\newcommand{\bP}{\mathds{P}}
\newcommand{\bR}{\mathds{R}}

\newcommand{\rD}{\mathscr{D}}

\newcommand{\sC}{\mathcal{C}}
\newcommand{\sD}{\mathcal{D}}
\newcommand{\sF}{\mathcal{F}}
\newcommand{\sH}{\mathcal{H}}
\newcommand{\sL}{\mathcal{L}}
\newcommand{\sU}{\mathcal{U}}

\newcommand{\mA}{\mathrm{A}}
\newcommand{\mC}{\mathrm{C}}
\newcommand{\mD}{\mathrm{D}}

\newcommand{\ud}{\mathrm{d}}
\newcommand{\udt}{\,\mathrm{d}}

\newcommand{\ind}{\mathbh{1}}
\newcommand{\aC}{\mathop{\sC}\limits^{\curvearrowleft}}

\renewcommand{\epsilon}{\varepsilon}
\renewcommand{\phi}{\varphi}
\makeatother

\begin{document}
\begin{frontmatter}

\title{An infinite-dimensional approach to path-dependent Kolmogorov equations}
\runtitle{Path-dependent Kolmogorov equations}

\begin{aug}
\author[A]{\fnms{Franco}~\snm{Flandoli}\ead[label=e1]{flandoli@dma.unipi.it}}
\and
\author[A]{\fnms{Giovanni}~\snm{Zanco}\corref{}\ead[label=e2]{zanco@mail.dm.unipi.it}}
\runauthor{F. Flandoli and G. Zanco}
\affiliation{Universit\`a di Pisa}
\address[A]{Dipartimento di Matematica\\
Universit\`a di Pisa\\
Largo Bruno Pontecorvo 5\\
56127 Pisa IT \\
Italy\\
\printead{e1}\\
\phantom{E-mail: }\printead*{e2}}
\end{aug}

%
\received{\smonth{12} \syear{2013}}
%
\revised{\smonth{4} \syear{2015}}

%
\begin{abstract}
In this paper, a Banach space framework is introduced in order to deal
with finite-dimensional path-dependent stochastic differential equations. A~version of
Kolmogorov backward equation is formulated and solved both in the space
of $L^p$ paths and in the space of continuous paths using the
associated stochastic differential equation, thus establishing a
relation between path-dependent SDEs and PDEs in analogy with the
classical case. Finally, it is shown how to establish a connection
between such Kolmogorov equation and the analogue
finite-dimensional equation that can be formulated in terms of the path-dependent
derivatives recently introduced by Dupire, Cont and Fourni\'e.
\end{abstract}

%
\end{frontmatter}

\section{Introduction}\label{secintro}
In the recent literature, a growing interest for path-depen\-dent
stochastic equations has arisen, due both to their mathematical
interest and to their possible applications in finance.

The path-dependent SDEs considered here will be of the form
%
\begin{equation}
\label{eqintroPSDE} \cases{ \ud X(t) = b_t(X_t)\udt t+\sigma\udt
W(t), &\quad for $t\in[t_0,T]$,
\cr
X_{t_0}=
\gamma_{t_0},}
\end{equation}
where $ \{W(t) \}_{t\geq0}$ is a Brownian motion in $\bR
^d$, $\sigma$ is a diagonalizable $d\times d$ matrix, the solution
$X(t)$ at time $t$ takes values in $\bR^d$, the notation $X_t$ stands
for the path of the solution on the interval $[0,t]$, $b_t$ is,  for
each $t\in[0,T]$, a map from a suitable space of paths to $\bR^d$,
$\gamma_{t_0}$ is a given path on $[0,t_0]$.

After the insightful ideas proposed by \citet{D} and
\citeauthor{CF1} (\citeyear{CF2,CF3,CF1}), who introduced a new concept of derivative and
developed a path-de\-pen\-dent It\^o formula which exhibits a first
connection between SDEs and PDEs in the path-dependent situation, some
effort was made into generalizing some classical concept to this
setting, like forward--backward systems and viscosity solutions [see
\citet{PW},
\citet{TZ},
\citet{EKTZ},
\citeauthor{ETZ1} (\citeyear{ETZ1,ETZ2}),
\citet{C1}]. Also, depending on the approach, there
are some similarities with investigations about delay equations; see,
for instance, \citet{FGG}, \citet{GM}, \citet{FMT}.
Some of these works formulate a path-dependent Kolmogorov equation
associated to the path-dependent SDE (\ref{eqintroPSDE}). Several
issues about such Kolmogorov equation are of interest. The purpose of
our work is to prove \emph{existence of classical $C^2$ solutions} and
to develop a Banach space framework suitable for this problem. To this
aim, we follow the classical method based on the probabilistic
representation formula in terms of solutions to the SDE, which however,
as explained in detail below, requires a new nontrivial analysis in our
framework.
\subsection{Notation}
We will use the following notation throughout the paper (in addition to
those introduced above): $T$ will stand for a fixed finite
time-horizon; $X_t(r)$ will stand again for the value of $X$ at $r$,
$r\leq t$. Stochastic processes will be denoted with upper-case
letters, while Greek lower-case letters will be used for deterministic
paths, most of the times seen as points in some paths space. As long as
no stochastics are involved, one can always think of a path $\gamma$
as defined on the whole interval $[0,T]$ and read $\gamma_t$ as its
restriction to $[0,t]$.

By $C([a,b];\bR^d)$ and $D([a,b];\bR^d)$ we will denote, respectively,
the space of continuous and c\`adl\`ag functions from the real interval
$[a,b]$ into $\bR^d$; $D([a,b);\bR^d)$ will denote the set of c\`adl\`
ag functions that have finite left limit also for $t\to b$.

\subsection{Main results}
\label{subsecmainresults}
A \emph{path-dependent nonanticipative function} is a family of
functions $b= \{b_t \}_{t\in[0,T]}$, each one being defined
on $D([0,t];\bR^d)$ with values in $\bR^d$ and being measurable with
respect to the canonical $\sigma$-field on $D([0,t];\bR^d)$.
Some possible examples of path-dependent functions are the following:
\begin{longlist}[(iii)]
\item[(i)] for $g\colon[0,T]\times[0,T]\times\bR^d\times
\bR^d\to\bR^d$ smooth, consider the function
\[
b_t (\gamma_t )=\int_0^tg
\bigl(t,s,\gamma(t),\gamma(s) \bigr)\udt s;
\]
\item[(ii)] for $0=t_0\leq t_1\leq t_2\leq\cdots\leq t_n\leq T$
fixed consider the function
\[
b_t (\gamma_t )=h_{i(t)} \bigl(\gamma(t),
\gamma(t_1),\dots,\gamma(t_{i(t)} ) \bigr),
\]
where for each $t\in[0,T]$ the index $i(t)\in\{0,\dots,n\}$ is such
that $t_{i(t)}\leq t<t_{i(t)+1}$ and, for each $j\in\{0,\dots,n\}$,
$h_j\colon\bR^{d\times(j+1)}\to\bR^d$ is a given function with
suitable properties;
\item[(iii)] for $\delta\in(0,T)$ and $q\colon\bR^{2d}\to
\bR^d$ smooth, consider the function
\[
b_t (\gamma_t )=q \bigl(\gamma(t),\gamma(t-\delta)
\bigr);
\]

\item[(iv)] in dimension $d=1$ consider the function
\[
b_t (\gamma_t )=\sup_{s\in[0,t]}\gamma(s).
\]
\end{longlist}

In order to formulate the path-dependent SDE (\ref{eqintroPSDE}) as an
SDE in Banach spaces, we consider it as a couple (endpoint, path) in
some infinite-dimensional space, as it is usually done for delay equations,
and reformulate consequently equation (\ref{eqintroPSDE}) as the
infinite-dimensional abstract SDE
%
\begin{equation}
\label{eqSDEintro}
\qquad\ud Y(t)=AY(t)\udt t+B\bigl(t,Y(t)\bigr)\udt t+\Sigma\,
\ud\beta(t)
\qquad\mbox{for }t\in[t_0,T],  Y(t_0)=y
\end{equation}
(understood in mild sense) where $A$ is the derivative operator, $B$ is
a sufficiently smooth (in Fr\'echet sense) nonlinear operator with
range in $\bR^d\times\{0\}$ and $\beta$ is a finite-dimensional
Brownian motion (Section~\ref{secstoch}).

We associate to it the backward Kolmogorov equation in integral form
with terminal condition $\Phi$
\begin{eqnarray}
\label{eqPKolmogorovintro} u ( t,y ) -\Phi( y ) &=& \int_{t}^{T}
\bigl\langle Du ( s,y ), Ay+B ( s,y ) \bigr\rangle\udt s
\nonumber\\[-8pt]\\[-8pt]\nonumber
&&{} +\frac{1}%
{2}\int_{t}^{T}\sum
_{j=1}^d\sigma_j^2
D^2u(s,y) (e_j,e_j)\udt s
\end{eqnarray}
and the related concept of solution (Section~\ref{seckolm}).

Our main result, under suitable regularity assumptions on $B$ and $\Phi
$, as explained in Section~\ref{secC} is the following (see Theorem
\ref{thmmain} for the precise statement):
\begin{theorem*}
The function
\[
u(s,y)=\bE\bigl[\Phi\bigl(Y^{s,y}(T) \bigr) \bigr],
\]
where $Y^{s,y}(t)$ solves equation (\ref{eqSDEintro}), is of class
$C^{2}$ with respect to $y$ and solves the backward Kolmogorov equation.
\end{theorem*}
Since under our assumptions all the integrands appearing in (\ref
{eqPKolmogorovintro}) are in $L^\infty$, a posteriori the function $u$
is Lipschitz in $t$ and hence, by Rademacher's theorem, differentiable
almost everywhere with respect to $t$. Therefore, for almost every $t$
it satisfies Kolmogorov backward equation in its differential form:
\[
\cases{
\displaystyle\frac{\partial u}{\partial t}(t,y)+\bigl\langle
Du(t,y),Ay+B(t,y)\bigr
\rangle+\frac{1}{2}\sum_{j=1}^d
\sigma^2_jD^2u(t,y) (e_j,e_j)=0,
\cr
\displaystyle u(T,\cdot)=\Phi.} %
\]
We moreover show (Section~\ref{secex}) that all usual examples
satisfy the regularity requirements of the previous theorem. Finally,
we provide some links between our results and the path-dependent
calculus developed by Cont and Fourni\'e (Section~\ref
{seccomparison}). In doing so, the main result we have (again under
some regularity assumptions compatible with those of the previous
theorem) is the following:
\begin{theorem*}
The function
\[
\nu_s(\gamma_s)=\bE\bigl[f \bigl(X^{\gamma_t}(T)
\bigr) \bigr],
\]
where $X^{\gamma_s}(t)$ is the solution to equation (\ref
{eqintroPSDE}), solves the path-dependent backward Kolmogorov equation
%
\begin{equation}
\label{eqPKolmogorovintroCF} %
\cases{
\displaystyle\rD_t\nu(\gamma_t)+b_t(
\gamma_t)\cdot\rD\nu_t(\gamma_t)+
\frac{1}{2}\sum_{j=1}^d
\sigma_j^2\rD^2_i
\nu_t(\gamma_t)=0,
\cr
\nu_T(
\gamma_T)=f(\gamma_T),} %
\end{equation}
in which the derivatives are understood as horizontal and vertical
derivatives as defined by \citet{CF1}.
\end{theorem*}
%
\subsection{Some ideas about the proofs}
We intend here to find regular solutions to the Kolmogorov\vspace*{1pt} equation, by
analogy with the classical theory. To this aim, the space of $L^2$
paths would appear to be the easiest setting to work in; unfortunately
there are no significant example of path-dependent functions, not even
integral functions, that satisfy the natural condition of having
uniformly continuous second Fr\'echet derivative in $L^2$; this is
discussed in detail in Section~\ref{secex}. To include a wider class
of functions, one would want to formulate and solve equations (\ref
{eqSDEintro}) and (\ref{eqPKolmogorovintro}) in the space of
continuous paths, that in our framework is the space
\[
\aC:= \Bigl\{y=\pmatrix{x
\cr
\phi}\in\bR^d\times C \bigl([-T,0
);\bR^d \bigr)\mbox{ s.t. } x=\lim_{s\uparrow0}\phi(s)
\Bigr\}.
\]
This leads to two issues: first, the operator $B$ (our abstract
realization of the functional $b$) takes values in a space larger than
$\aC$, thus we have to consider paths with a single jump-discontinuity
at the final time $t=0$. But then the semigroup generated by $A$ shifts
such discontinuity so that we have to deal with paths with a single
discontinuity at an arbitrary time $t$. The need to work with a linear
space and possibly with a Banach space structure suggests the choice of
\[
\sD:=\bR^d\times D \bigl([-T,0);\bR^d \bigr)
\]
with the uniform norm as the ambient space for our equations.

The second issue comes along when we try to establish the link between
the SDE and the PDE. As in the classical theory, we need to work with
some It\^o-type formula. We decide not to use some version of the It\^o
formula in Banach spaces due to the difficulties one encounters in
defining a concept of quadratic variation there [see, e.g.,
\citeauthor{DRth} (\citeyear{DRth,DR2,DR1}), \citet{DFR}], although we intend to address
this problem in our future works; we proceed therefore using a Taylor
expansion, but we are not able to control the second-order terms in
spaces endowed with the uniform norm.

Therefore, we adopt the following strategy: we go back to an $L^p$
setting with $p\geq2$ (recovering in this way at least examples like
integral functionals) and we develop rigorously the relation between
the SDE and the PDE in this framework (Section~\ref{secLp}). We then
introduce an approximation procedure to extend our results to the space
of continuous paths (Section~\ref{secC}). This step requires us to
introduce an additional assumption that remarks again the deep effort
that is needed in order to obtain a satisfactory general theory already
in the easiest case of regular coefficients.

\section{The stochastic equation}
\subsection{Framework}
\label{secstoch}
From now onward, fix a time horizon $0<T<\infty$ and a filtered
probability space $ (\Omega,\sF, \{\sF_t \}_{t\in
[0,T]},\bP)$.
We introduce the following spaces:
\begin{eqnarray*}
\sC&:=&\bR^d\times\Bigl\{\phi\in C_b \bigl([-T,0
);\bR^d \bigr) \colon\exists\lim_{s\uparrow0}\phi(s) \Bigr
\},
\\
\aC&:=& \Bigl\{y=\pmatrix{x
\cr
\phi}\in\sC\mbox{ s.t. }x=\lim
_{s\uparrow0}\phi(s) \Bigr\},
\\
\sD&:=&\bR^d\times D_b \bigl([-T,0);
\bR^d \bigr),
\\
\sD_t&:=&\left\{y=\pmatrix{x
\cr
\phi}\in\sD\mbox{ s.t. }\phi\mbox{ is
discontinuous at most in the only point $t$}\right\},
\\
\sL^p&:=&\bR^d\times L^p \bigl(-T,0;
\bR^d \bigr),\qquad p\geq2.
\end{eqnarray*}
All of them apart from $\sL^p$ are Banach spaces with respect to the
norm $\llVert {x \choose \phi}\rrVert ^2=\llvert x\rrvert
^2+{\llVert \phi\rrVert }_\infty
^2$, while $\sL^p$ is a Banach space with respect to the norm $\llVert
{x \choose  \phi}\rrVert ^2=\llvert x\rrvert ^2+\llVert \phi
\rrVert
_{p}^2$; the space
$\sD$ turns out to be not separable with respect to this norm but this
will not undermine our method.

With these norms, we have the natural relations
\[
\aC\subset\sC\subset\sD\subset\sL^p
\]
with continuous embeddings. We remark that $\aC$, $\sC$ and $\sD$
are dense in $\sL^p$ while neither $\aC$ nor $\sC$ are dense in $\sD
$. The choice for the interval $[-T,0]$ is made in accordance with most
of the classical literature on delay equations.

Notice that the space $\aC$ has not the structure of a product space;
notice also that it is isomorphic to the space $C ([-T,0];\bR
^d )$.

As said above, we consider a family $b= \{b_t \}_{t\in
[0,T]}$ of functions
\[
b_t:D \bigl([0,t];\bR^d \bigr)\rightarrow
\bR^d
\]
adapted to the canonical filtration and we formulate the path-dependent
stochastic differential equation
%
\begin{equation}
\label{eqPSDE} \cases{ \ud X(t) = b_t(X_t)\udt t+\sigma\udt
W(t), &\quad for $t\in[t_0,T]$,
\cr
X_{t_0} =
\gamma_{t_0},}
\end{equation}
where $\sigma$ is a diagonalizable $d\times d$ matrix and $W$ is a
$d$-dimensional Brownian motion. $b$ can also be seen as an $\bR
^d$-valued function on the space $D=\bigcup_t D ([0,t];\bR^d
)$.

To reformulate the path-dependent SDE (\ref{eqPSDE}) in our framework,
we need to introduce two linear bounded operators: for every $t\in
[0,T]$ define the \emph{restriction operator}
%
\begin{eqnarray}
\label{eqMt}
\nonumber
&\displaystyle M_t\colon D \bigl([-T,0);
\bR^d \bigr)\longrightarrow D \bigl([0,t);\bR^d \bigr),&
\nonumber\\[-8pt]\\[-8pt]\nonumber
&\displaystyle M_t(\phi) (s)=\phi(s-t),\qquad{s\in[0,t)}, &
\end{eqnarray}
and the \emph{backward extension} operator
%
\begin{eqnarray}
\label{eqLt}
\nonumber
& \displaystyle L_t\colon D \bigl([0,t);
\bR^d \bigr)\longrightarrow D \bigl([-T,0);\bR^d \bigr),&
\nonumber\\[-8pt]\\[-8pt]\nonumber
&\displaystyle L_t(\gamma) (s)=\gamma(0)\ind_{[-T,-t)}(s)+\gamma(t+s)\ind
_{[-t,0)}(s),\qquad{s\in[-T,0)}.&
\end{eqnarray}
Since the extension in the definition of $L_t$ is arbitrary, one has that
%
\begin{equation}
\label{eqML} M_tL_t\gamma=\gamma
\end{equation}
while in general
\[
L_tM_t\phi\neq\phi.
\]
Note also that both $L_t$ and $M_t$ map continuous functions into
continuous functions.
Set moreover
%
\begin{equation}
\label{eqtildeMt} \widetilde{M}_t\pmatrix{x
\cr
\phi}(s)= \cases{
M_t\phi(s), &\quad$s\in[0,t)$,
\cr
x, &\quad$s=t$.}
\end{equation}
Now given a functional $b$ as in (\ref{eqPSDE}) one can define a
function $\hat b$ on $[0,T]\times\sD$ setting
%
\begin{equation}
\label{eqB} \hat b \left(t,\pmatrix{x
\cr
\phi} \right)=\hat b(t,x,\phi):=
b_t \left(\widetilde{M}_t\pmatrix{x
\cr
\phi} \right);
\end{equation}
conversely if $\hat b$ is given one can obtain a functional $b$ on $D$ setting
%
\begin{equation}
\label{eqb} b_t(\gamma): =\hat b\bigl(t,\gamma(t),L_t
\gamma\bigr).
\end{equation}
The idea is simply to shift and extend (or restrict) the path in order
to pass from one formulation to another.

For instance, the functional of example (i) in Section~\ref{secintro}
would define a function $\hat b$ on $[0,T]\times\sD$
given by
%
\begin{equation}
\label{ex1bis} \hat b \biggl(t,\pmatrix{x
\cr
\phi} \biggr)=\int_0^tg
\bigl(t,s,x,\phi(s-t) \bigr)\udt s.
\end{equation}
We consider again the path-dependent SDE (\ref{eqPSDE}) with the
initial condition given now by a path $\psi$ on $[-T+t_0,t_0]$ and its
terminal value $x=\psi(t_0)$,
%
\begin{equation}
\label{eqPSDE2} \cases{ \ud X(s) = b_s(X_s)\udt s+
\sigma\udt  W(s), &\quad for $s\in[t_0,T]$,
\cr
X(t_0) = x=
\psi(t_0),
\cr
X(s) = \psi(s), &\quad for $s\in[-T+t_0,t_0)$.}
\end{equation}
Recall that by $X_s$ we denote the path of $X$ starting from $0$ up to
time $s$, not a portion of the path of $X$ of length $T$, which would
be anyway well defined in this setting.
If $X$ solves (\ref{eqPSDE2}) (in some space), for $t\in[t_0,T]$ we set
\[
Y(t)=\pmatrix{ X(t)
\cr
\bigl\{X(t+s) \bigr\}_{s\in[-T,0]}}
\]
and then differentiate with respect to $t$ formally obtaining
%
\begin{equation}
\label{eqDY} \qquad\frac{\ud Y(t)}{\ud t}=\pmatrix{ \dot{X}(t)
\cr
\bigl\{\dot{X}(t+s)
\bigr\}_s} =\pmatrix{ 0
\cr
\bigl\{\dot{X}(t+s) \bigr
\}_s}+ \pmatrix{b_t(X_t)
\cr
0}+\pmatrix{
\sigma\dot{W}(t)
\cr
0}.
\end{equation}
It is therefore natural to define the operators
%
\begin{eqnarray}
\label{eqdefA} A\pmatrix{x
\cr
\phi}&: =& \pmatrix{0
\cr
\dot\phi},
\\
\label{eqdefB} B\biggl(t,\pmatrix{x
\cr
\phi}\biggr)&: =& \pmatrix{\hat b \biggl(t,\pmatrix{x
\cr
\phi} \biggr)
\vspace*{3pt}\cr
0}
\end{eqnarray}
and
%
\begin{equation}
\label{eqdefSigma} \Sigma\pmatrix{x
\cr
\phi}: = \pmatrix{\sigma x
\cr
0}
\end{equation}
and to formulate the infinite-dimensional SDE
%
{\renewcommand{\theequation}{14$^\prime$}
\begin{equation}
\label{eqDYshort} \ud Y(t)=AY(t)\udt t+B \bigl(t,Y(t) \bigr)\udt t+\Sigma\udt\beta(t), \qquad t\in[t_0,T], 
\end{equation}\setcounter{equation}{17}}%
where $\beta$ is given by
%
\begin{equation}
\label{eqdefbeta} \beta(t)=\pmatrix{W(t)
\cr
0},
\end{equation}
with some initial condition $Y(t_0)=y$.

Solutions of this SDE will always be understood to be \emph{mild
solutions}, that is, we want to solve
%
{\renewcommand{\theequation}{14$^{\prime\prime}$}
\begin{equation}
\label{eqDYmild} Y(t)=e^{(t-t_0)A}y+\int_{t_0}^te^{(t-s)A}B
\bigl(s,Y(s) \bigr)\udt s+\int_{t_0}^te^{(t-s)A}
\Sigma\udt\beta(s).
\end{equation}\setcounter{equation}{18}}%
It is not difficult to show that if $Y$ solves (\ref{eqDYshort})
then its first coordinate $X(t)$ solves the original SDE (\ref
{eqPSDE2}).

\subsection{Some properties of the convolution integrals}
\label{secproperties}
The operator $A$ has different domains depending on the space that we work in; we set
\begin{eqnarray*}
\Dom(A) &=& \biggl\{\pmatrix{x
\cr
\phi}\in\sL^p:\phi\in
W^{1,p} \bigl( -T,0;\bR^d \bigr), \phi( 0 )=x \biggr\},
\\
\Dom(A_{\aC}) &=& \biggl\{\pmatrix{x
\cr
\phi}\in\aC:\phi\in
C^1 \bigl([-T,0);\bR^d \bigr) \biggr\};
\end{eqnarray*}
one can think to define $A$ on $\sL^p$ and then consider its
restriction to $\sD$ or to $\aC$, as the notation above emphasizes.

It is well known [see Theorem 4.4.2 in \citet{BDPDM1}] that $A$
is the infinitesimal generator of a strongly continuous semigroup both
in $\sL^p$ and in $\aC$; it is easy to check that it still generates
a semigroup in $\sD$ which is not uniformly continuous. Indeed we have that
%
\begin{equation}
\label{eqetA} e^{tA}\pmatrix{ x
\cr
\phi} = \pmatrix{ x
\cr
\bigl\{\phi
( \xi+t ) \ind_{ [ -T,-t ) } ( \xi)+x\ind_{ [ -t,0 ] } ( \xi) \bigr
\}_{\xi\in[-T,0]}}.
\end{equation}
This formula comes from the trivial delay equation
\[
\cases{ \displaystyle\frac{\mathrm{d}x ( t ) }{\ud t} =0, &\quad$t\geq0$,
\vspace*{5pt}\cr
x ( 0 ) =x, & \quad$x
( \xi) =\varphi(\xi)$ for $\xi\in[ -T,0 ]$;}
\]
its solution, for $t\geq0$, is simply $x ( t ) =x$. If we introduce
the pair
\[
y ( t ):= \pmatrix{ x ( t )
\cr
x_{\mid_{ [ t-T,t ] }}}
\]
then
\[
y ( t ) =e^{tA} \pmatrix{ x
\cr
\varphi}.
\]
However, it still holds that
%
\begin{equation}
\label{eqnormetA} \bigl\llVert e^{tA}\bigr\rrVert_{L(\sD,\sD)}\leq C
\qquad\mbox{for }t\in[0,T]
\end{equation}
with $C$ not depending on $t$. Moreover, it is evident from (\ref
{eqetA}) that $e^{tA}$ maps $\sL^p$ into $\sL^p$, $\sD$ into $\sD$
and $\aC$ into $\aC$, but it maps $\sC$ into $\sD_{-t}$ because an
element of $\sC$ is essentially a continuous function with a unique
discontinuity at its endpoint, and the semigroup just shifts that
discontinuity. In particular this happens for elements of $\bR^d\times
\{0\}$.

Consider the stochastic convolution
\[
Z^{t_{0}} ( t ):=\int_{t_{0}}^{t}e^{ ( t-s )
A}
\Sigma\udt\beta( s ) =\int_{t_{0}}^{t}e^{ (
t-s )
A}
\pmatrix{\sigma\udt  W ( s )
\cr
0}, \qquad t\geq t_{0}.
\]
It is not obvious to investigate $Z^{t_0}$ by infinite-dimensional
stochastic integration theory, due to the difficult nature of the
Banach space $\sD$. However, we may study its properties thanks to the
following explicit formulas.
From now on, we work in a set $\Omega_0\subseteq\Omega$ of full
probability on which $W$ has continuous trajectories. For $\omega\in
\Omega_0$ fixed, for any $x\in\bR^d$ we have
\[
e^{ ( t-s ) A}\Sigma\pmatrix{x
\cr
0}= \pmatrix{ \sigma x
\cr
\bigl\{\sigma x
\ind_{ [ - ( t-s ), 0 ] } ( \xi) \bigr\}_{\xi
\in[ -T,0 ]}}
\]
hence
%
\begin{eqnarray}
\label{eqZ^t0} Z^{t_{0}} ( t ) &=& \pmatrix{ \displaystyle\int_{t_{0}}^{t}
\sigma\udt  W ( s )
\cr
\displaystyle\int_{t_{0}}^{t}
\ind_{ [ - ( t-s ), 0 ] } ( \cdot) \sigma\udt  W ( s )}
\nonumber\\[-8pt]\\[-8pt]\nonumber
\nonumber
&=& \pmatrix{ \displaystyle\sigma\bigl(W ( t ) -W ( t_{0} ) \bigr)
\cr
\displaystyle\sigma\bigl(W \bigl( ( t+\cdot) \vee t_{0} \bigr) -W (
t_{0} ) \bigr)}
\end{eqnarray}
because
\[
\int_{t_{0}}^{t}\ind_{ [ - ( t-s ), 0 ] } (\xi) \sigma
\udt W ( s ) =\int_{t_{0}}^{t}\ind_{ [ 0,t+\xi]} ( s
) \sigma\udt  W ( s ).
\]
From the previous formula, we see that $Z^{t_{0}} ( t )\in
\aC$, hence $Z^{t_0}(t)\in\sL^p$.

We have
\[
\bigl\llVert Z^{t_{0}} ( t ) \bigr\rrVert_{\aC}=2 \sup
_{\xi\in[ -T,0 ] }\bigl\llvert\sigma\bigl(W \bigl( ( t+\xi) \vee
t_{0} \bigr) -W ( t_{0} ) \bigr)\bigr\rrvert
\]
hence [using the fact that $r\mapsto W ( t_{0}+r ) -W
(t_{0} ) $ is a Brownian motion and applying Doob's inequality]
%
\begin{eqnarray}
\label{eqEZ} \mathbb{E} \bigl[ \bigl\llVert Z^{t_{0}} ( t ) \bigr
\rrVert_{\aC}^{4} \bigr] &\leq&2^4\mathbb{E}
\Bigl[ \sup_{s\in
[0,t-t_{0} ] }\bigl\llvert\sigma W ( s ) \bigr\rrvert
^{4} \Bigr]
\nonumber\\[-8pt]\\[-8pt]\nonumber
&\leq& C^{\prime}\mathbb{E} \bigl[ \bigl\llvert W (
t-t_{0} ) \bigr\rrvert^{4} \bigr] \leq C^{\prime
\prime} (
t-t_{0} ) ^{2},
\end{eqnarray}
where $C^{\prime}$ and $C^{\prime\prime}$ are suitable constants.
Consequently, the same property holds in $\sL^p$ (possibly with a
different constant) by continuity of the embedding $\aC\subset\sL^p$.
Moreover, from (\ref{eqZ^t0}) we obtain that for $\omega$ fixed
\begin{eqnarray*}
&& \bigl\llVert Z^{t_0}(t)-Z^{t_0}(s)\bigr\rrVert
_{\aC}
\\
&&\qquad =C \Bigl(\bigl\llvert W(t)-W(s)\bigr\rrvert+\sup_{\xi\in[-T,0]}\bigl
\llvert W \bigl((t+\xi)\vee t_0 \bigr)-W \bigl((s+\xi)\vee
t_0 \bigr)\bigr\rrvert\Bigr).
\end{eqnarray*}
Observe that (we suppose $s<t$ for simplicity)
\begin{eqnarray*}
&& W \bigl((t+\xi)\vee t_0 \bigr)-W \bigl((s+\xi)\vee t_0
\bigr)
\\
&&\qquad = \cases{ 0, &\quad$\xi\in[-T,t_0-t]$,
\cr
W(t+
\xi)-W(t_0), &\quad$\xi\in[t_0-t,t_0-s]$,
\cr
W(t+\xi)-W(s+\xi), &\quad$\xi\in[t_0-s,0]$} %
\end{eqnarray*}
and
\[
\sup_{\xi\in[t_0-t,t_0-s]}\bigl\llvert W(t+\xi)-W(t_0)\bigr
\rrvert=\sup_{\eta
\in[t_0,t_0+(t-s)]}\bigl\llvert W(\eta)-W(t_0)
\bigr\rrvert,
\]
therefore, $Z^{t_0}$ is a continuous process in $\aC$, since any fixed
trajectory of $W$ is uniformly continuous. The same property holds then
in $\sL^p$ again by continuity of the embedding $\aC\subset\sL^p$.
We can argue in a similar way for $F^{t_0}:[t_0,T]\times L^\infty
([t_0,T];\sD)\rightarrow\sD$,
\[
F^{t_{0}} (t,\theta) =\int_{t_{0}}^{t}e^{ (
t-s )A}B
\bigl( s,\theta( s ) \bigr)\udt s.
\]
From (\ref{eqdefB}), using (\ref{eqetA}) one deduces that
\[
e^{(t-s)A}B \bigl(s,\theta(s) \bigr)= \pmatrix{b_s \bigl(
\widetilde M_s\theta(s) \bigr)
\vspace*{3pt}\cr
b_s \bigl(\widetilde
M_s\theta(s) \bigr)\ind_{[-t+s]}(\xi)}
\]
and, therefore,
\[
\int_{t_0}^te^{(t-s)A}B \bigl(s,\theta(s)
\bigr)\udt s= \pmatrix{ \displaystyle\int_{t_0}^t
b_s \bigl(\widetilde M_s\theta(s) \bigr)\udt s
\vspace*{3pt}\cr
\displaystyle\biggl\{\int_{t_0}^{t+\xi}b_s \bigl(
\widetilde M_s\theta(s) \bigr)\udt s \biggr\}_\xi}
\]
which shows that $F^{t_0}(t,\theta)$ always belongs to $\aC$. Writing
\[
Y^{t_0,y}(t)=e^{(t-t_0)A}y+F^{t_0}\bigl(t,Y^{t_0,y}
\bigr)+Z^{t_0}(t)
\]
we see immediately that, for any $t\in[t_0,T]$, $Y^{t_0,y}(t)\in\sD$
if $y\in\sD$ and\break $Y^{t_0,y}(t)\in\aC$ if $y\in\aC$. This will be
crucial in the sequel.

\subsection{Existence, uniqueness and differentiability of solutions
to the SDE}
We state and prove here some abstract results about existence and
differentiability of solutions to the stochastic equation
%
{\renewcommand{\theequation}{14$^\prime$}
\begin{equation}
\ud Y(t)=AY(t)\udt t+B \bigl(t,Y(t) \bigr)\udt t+\Sigma\udt\beta(t),
\qquad Y(t_0)=y,
\end{equation}\setcounter{equation}{22}}%
with respect to the initial data. By abstract, we mean that we consider
a general $B$ not necessarily defined through a given $b$ as in
previous sections. Also $A$ can be thought here to be a generic
infinitesimal generator of a semigroup which is strongly continuous in
$\sL^p$ and satisfies (\ref{eqnormetA}) in $\sD$. Although all
these theorems are analogous to well-known results for stochastic
equations in Hilbert spaces [see, e.g., \citet{DPZrosso}], we
give here complete and exact proofs due to the lack of them in the
literature for the case of time-dependent coefficients in Banach
spaces, which is the one of interest here.

We are interested in solving the SDE in $\sL^p$ and in $\sD$; since
almost all the proofs can be carried out in the same way for each of
the spaces we consider and since we do not need any particular property
of these spaces themselves, we state all our results in this section in
a general Banach space $E$, stressing out possible distinctions that
could arise from different choices of $E$. In the following, we will
identify $L (E,L(E,E) )$ with $L(E,E;E)$ (the space of bilinear
forms on $E$) in the usual way.

We will make the following assumption.
\begin{assumption}
\label{assB}
\[
B\in L^{\infty} \bigl( 0,T;C_{b}^{2,\alpha} (E,E ) \bigr)
\]
for some $\alpha\in( 0,1 ) $, where we have denoted by
$C_{b}^{2,\alpha} (E,E ) $ the space of twice Fr\'{e}chet
differentiable functions $\varphi$ from $E$ to $E$, bounded with their
differentials of first and second order, such that $x\mapsto
D^{2}\varphi( x ) $ is $\alpha$-H\"{o}lder continuous
from $E$ to $L (E,E;E ) $. The $L^{\infty}$ property in
time means that the differentials are measurable in $ ( t,x
) $ and both the function, the two differentials and the H\"{o}lder
norms are bounded in time. Under these conditions, $B$, $DB$, $D^{2}B$
are globally uniformly continuous on $E$ [with values in $E$, $L
(E,E ) $, $L (E,E;E ) $], respectively, and with a
uniform in time modulus of continuity.
\end{assumption}
\begin{theorem}
\label{thmSDE}
Equation (\ref{eqDYshort}) can be solved in a mild sense path by
path: for any $y\in E$, any $t_0\in[0,T]$ and every $\omega\in\Omega
_0$ there exists a unique function $[t_0,T]\ni t\to Y^{t_0,y}(t,\omega
)\in E$ which satisfies identity (\ref{eqDYmild})
{\renewcommand{\theequation}{14$^{\prime\prime}$}
\begin{eqnarray}
Y^{t_0,y}(t,\omega)&=&e^{(t-t_0)A}y+\int
_{t_0}^te^{(t-s)A}B \bigl(s,Y^{t_0,y}(s,
\omega) \bigr)\udt s
\nonumber\\[-8pt]\\[-8pt]\nonumber
&&{}+\int_{t_0}^te^{(t-s)A}\Sigma\udt\beta(s,
\omega).
\end{eqnarray}\setcounter{equation}{22}}%
Such a function is continuous if $E=\sL^p$, it is only in $L^\infty$
if $E=\sD$.
\end{theorem}

\begin{pf}
Thanks to the Lipschitz property of $B$ the proof follows through a
standard argument based on the contraction mapping principle. The lack
of continuity in $\sD$ is due to the fact that the semigroup $e^{tA}$
is not strongly continuous in~$\sD$.
\end{pf}
\begin{theorem}
\label{4535454754545}
For every $\omega\in\Omega_0$, for all $t_0\in[0,T]$ and $t\in
[t_0,T]$ the map $y\mapsto Y^{t_0,y}(t,\omega)$ is twice Fr\'echet
differentiable and the map\break $y\mapsto D^{2}Y^{t_{0},y} ( t,\omega
) $ is $\alpha$-H\"{o}lder continuous from $E$ to $L (
E,E;E )$. Moreover, if $E=\sL^p$, for any fixed $t$ and $y$ the
map $s\mapsto Y^{s,y}(t,\omega)$ is continuous. If $E=\sD$, the same
conclusion holds only for any fixed $y\in\aC$.
\end{theorem}
\begin{pf}
Due to its length the proof is postponed to the \hyperref[append]{Appendix}.
\end{pf}
\begin{theorem}
\label{thmmarkov}
If the solution $Y^{t_0,y}(t)$ is continuous as a function of $t$ with
values in $E$, then it has the Markov property.
\end{theorem}
\begin{pf}
This follows immediately from Theorem 9.15 on \citet{DPZrosso}.
Notice that there the authors require a different set of hypothesis
which, however, are needed only for proving existence and uniqueness of
solutions and not in the actual proof of the result. It therefore
applies to our situation as well.
\end{pf}
In Section~\ref{secLp}, we will need the notion of modulus of
continuity for the second Fr\'echet derivative of a map from $E$ into
$E$, together with some of its properties; we summarize what we will
need in the following general remark.
\begin{remark}
\label{RemTaylor}Given a map $R:E\rightarrow L ( E,E;\mathbb{R}%
) $, we define its modulus of continuity
\[
\mathfrak{w} ( R,r ) =\sup_{{\llVert y-y^{\prime
}\rrVert }_{E}\leq r}{\bigl\llVert R ( y ) -R
\bigl( y^{\prime} \bigr) \bigr\rrVert}_{L ( E,E;\mathbb{R} ) }.
\]
Let $v:E\rightarrow\mathbb{R}$ be a function with two Fr\'{e}chet
derivatives at each point, uniformly continuous on bounded sets. Then
there exists a function $r_{v}:E^{2}\rightarrow\mathbb{R}$ such that
\begin{eqnarray*}
v ( x ) -v ( x_{0} )  &=& \bigl\langle Dv ( x_{0} )
,x-x_{0} \bigr\rangle+\tfrac{1}{2}D^{2}v (
x_{0} ) ( x-x_{0},x-x_{0} ) +
\tfrac{1}{2}r_{v} ( x,x_{0} ),
\\
\bigl\llvert r_{v} ( x,x_{0} ) \bigr\rrvert& \leq&
\mathfrak{w} \bigl( D^{2}v,{\llVert x-x_{0}\rrVert
}_{E} \bigr){\llVert x-x_{0}\rrVert}_{E}^{2}
\end{eqnarray*}
for every $x,x_{0}\in E$. Indeed,
\[
v ( x ) -v ( x_{0} ) = \bigl\langle Dv ( x_{0} )
,x-x_{0} \bigr\rangle+\tfrac{1}{2}D^{2}v ( \xi
_{x,x_{0}} ) ( x-x_{0},x-x_{0} ),
\]
where $\xi_{v,x,x_{0}}$ is an intermediate point between $x_{0}$ and $x$,
and thus
\begin{eqnarray*}
\bigl\llvert r_{v} ( x,x_{0} ) \bigr\rrvert& =&\bigl
\llvert\bigl( D^{2}v ( \xi_{v,x,x_{0}} ) -D^{2}v
(x_{0} ) \bigr) ( x-x_{0},x-x_{0} ) \bigr
\rrvert
\\
& \leq& \bigl\llVert D^{2}v ( \xi_{v,x,x_{0}} )-D^{2}v
( x_{0} ) \bigr\rrVert_{L ( E,E;\mathbb{R} )
}{\llVert x-x_{0}
\rrVert}_{E}^{2}
\\
& \leq&\mathfrak{w} \bigl( D^{2}v,{\llVert x-x_{0}\rrVert
}_{E} \bigr){\llVert x-x_{0}\rrVert}_{E}^{2}
.
\end{eqnarray*}
If $D^{2}v$ is $\alpha$-H\"older continuous, namely
\[
{\bigl\llVert D^2v ( y ) -D^2v \bigl( y^{\prime}
\bigr) \bigr\rrVert}_{L ( E,E;\mathbb{R} ) }\leq M{\bigl\llVert
y-y^{\prime}\bigr
\rrVert}_{E}^{\alpha}%
\]
then
\[
\mathfrak{w} \bigl( D^{2}v,{\llVert x-x_{0}\rrVert
}_{E} \bigr) \leq M{\llVert x-x_{0}\rrVert
}_{E}^{\alpha}%
\]
and thus
\[
\bigl\llvert r_{v} ( x,x_{0} ) \bigr\rrvert\leq M{
\llVert x-x_{0}\rrVert}_{E}^{2+\alpha}.
\]
\end{remark}
%
\section{The Kolmogorov equation}
\label{seckolm}
In this and the following two sections, we introduce and solve the
backward Kolmogorov equation in our infinite-dimensional setting. The
relation between the results we shall show and the finite-dimensional
path-dependent SDE we started from will be investigated in Section~\ref
{seccomparison}.

Suppose for a moment we are working in a standard Hilbert-space
setting, that is, in a space $\sH=\bR^d\times H$ where $H$ is a
Hilbert space. Then [see again \citet{DPZrosso}] the backward
Kolmogorov equation, for the unknown $u:[0,T]\times\sH\to\bR$, is
%
\begin{equation}
\label{eqKolmstraight} %
\cases{ \displaystyle\frac{\partial u}{\partial t}(t,y)+\frac{1}{2}\tr
\bigl( \Sigma^\ast\Sigma D^{2}u(t,y) \bigr)+ \bigl\langle
Du(t,y),Ay+B ( t,y ) \bigr\rangle=0,
\vspace*{3pt}\cr
u(T,\cdot)=\Phi,} %
\end{equation}
where $\Phi$ is a given terminal condition and $Du$, $D^2u$ represent
the first and second Fr\'echet differentials with respect to the
variable $y$. Its solution, under suitable hypothesis on $A$, $B$,
$\Sigma$ and $\Phi$, is given by
%
\begin{equation}
\label{eqsolKolmHilbert} u(t,y)=\bE\bigl[\Phi\bigl(Y^{t,y}(T) \bigr)
\bigr],
\end{equation}
where $Y^{t,y}(t)$ solves the associated SDE
%
{\renewcommand{\theequation}{14$'$ bis}
\begin{eqnarray}
\label{eqSDEHilbert}
\ud Y(s)= \bigl[AY(s)+B \bigl(s,Y(s)
\bigr) \bigr]\udt s+\Sigma\udt\beta(s),
\nonumber\\[-8pt]\\[-8pt]
\eqntext{s\in[t,T],  Y(t)=y}
\end{eqnarray}\setcounter{equation}{24}}%
in $\sH$.
In our framework, where the spaces are only Banach spaces, we have to
give a precise meaning to the Kolmogorov equation and prove its
relation above with the SDE.

As outlined in the \hyperref[secintro]{Introduction}, we would like to solve it on the space
$\aC$, but since $B(t,y)$ belongs to $\bR^d\times\{0\}\nsubseteq\aC
$, in order to give meaning to the term $\langle Du(t,y),B(t,y)\rangle
$ we need $Du(t,y)$ to be a functional defined at least on $\sC$,
which necessarily implies $u$ to be defined on $[0,T]\times\sC$.
Therefore, we should solve (in mild sense) the SDE for $y\in\sC$ and
this implies that $Y^{t,y}(s)\in\sD_{-t+s}$ for $s\neq t$; this in
turn requires $\Phi$ to be defined at least on $\bigcup_{s\in[t,T]}\sD
_{-t+s}$ in order for a function of the form (\ref{eqsolKolmHilbert})
to be well defined. However, the space $\bigcup\sD_s$ is not a linear
space, thus it turns out that it is more convenient, also for
exploiting a Banach space structure, to formulate everything in $\sD$,
that is,
\[
u:[0,T]\times\sD\rightarrow\bR.
\]
Therefore, we interpret $\langle\cdot,\cdot\rangle$ in this setting
as the duality pairing between $\sD^\prime$ and~$\sD$.

For the trace term, if we denote by $e_1,\dots,e_d$ an orthonormal
basis of $\bR^d$ where $\sigma$ diagonalizes, that is, $\sigma
e_j=\sigma_j e_j$ for some real $\sigma_j$ (in any of the spaces
considered up to now), we could complete it to an orthonormal system $\{
e_n\}$ in $\sH$ obtaining that
\[
\tr\bigl(\Sigma^\ast\Sigma D^2u(t,y) \bigr)=\sum
_j\sigma_j^2\bigl
\langle D^2u(t,y)e_j,e_j\bigr\rangle;
\]
hence, by analogy, also when working in $\sD$ we interpret the trace
term as
%
\begin{equation}
\tr\bigl(\Sigma^\ast\Sigma D^2u(t,y) \bigr)=\sum
_{j=1}^d\sigma_j^2
D^2u(t,y) (e_j,e_j).
\end{equation}
Moreover, we consider Kolmogorov equation in its integrated form with
respect to time, that is, given a (sufficiently regular; see below)
real function $\Phi$ on $\sD$ we seek for a solution of the PDE:
\begin{eqnarray}
\label{eqPKolmogorov} u ( t,y ) -\Phi( y ) &=&\int_{t}^{T}
\bigl\langle Du ( s,y ), Ay+B ( s,y ) \bigr\rangle\udt s
\nonumber\\[-8pt]\\[-8pt]\nonumber
&&{}+\frac{1}{2}\int_{t}^{T}\sum
_{j=1}^d\sigma_j^2
D^2u(s,y) (e_j,e_j)\udt s.
\end{eqnarray}
Here, one can see one of the difficulties in working with Banach
spaces: the second-order term in the equation comes from the quadratic
variation of the solution of the SDE, but in such spaces there is no
general way of defining a quadratic variation [although, as mentioned
at the beginning, a general theory of quadratic variation is currently
being developed by F. Russo and collaborators; see the works
\citet{DR1}, \citet{DFR} and the references therein].

Although we will seek for such a $u$, when dealing with the equation we
will always choose $y$ to be in $\Dom(A_{\aC})$, to let all the terms
appearing there be well defined.

All these observations lead to our definition of solution to (\ref
{eqPKolmogorov}); first, we say that a functional $u$ on $[0,T]\times
\sD$ belongs to
\[
L^{\infty} \bigl( 0,T;C_{b}^{2,\alpha} (\sD,\bR) \bigr)
\]
if it is twice Fr\'echet differentiable on $\sD$, $u$, $Du$ and $D^2u$
are bounded, the map $x\mapsto D^2u(x)$ is $\alpha$-H\"older
continuous from $\sD$ to $L (\sD,\sD;\sD) $ (the space
of bilinear forms on $\sD$), the differentials are measurable in
$ ( t,x ) $ and the function, the two differentials and the
H\"{o}lder norms are bounded in time.
\begin{definition}
\label{defsolD}
Given $\Phi\in C_{b}^{2,\alpha} (\sD,\mathbb{R} ) $, we
say that $u: [ 0,T ] \times\sD\rightarrow\mathbb{R}$ is
a classical solution of the Kolmogorov equation with terminal condition
$\Phi$ if
\[
u\in L^{\infty} \bigl( 0,T;C_{b}^{2,\alpha} ( \sD,
\mathbb{R}%
) \bigr) \cap C \bigl( [ 0,T ] \times\aC,\mathbb
{R}%
\bigr),
\]
$u (\cdot,y )$ is Lipschitz for any $y\in\Dom
(A_{\aC} )$ and satisfies identity (\ref{eqPKolmogorov}) for
every $t\in[ 0,T ] $ and $y\in\Dom(A_{\aC}
) $, with the duality terms understood with respect to the topology of
$\sD$.
\end{definition}
It will be clear in Section~\ref{secC} that the restriction $y\in
\Dom(A_{\aC} )$ is necessary and that it would not be
possible to obtain the same result choosing $y$ in some larger space.

Our aim is to show that, in analogy with the classical case, the function
\[
u(t,y)=\bE\bigl[\Phi\bigl(Y^{t,y}(T) \bigr) \bigr]
\]
solves equation (\ref{eqPKolmogorov}).

However, we are not able to prove this result directly, due essentially
to the lack of an appropriate It\^o-type formula for our setting.
Therefore, we will proceed as follows: first, we are going to show how
to prove such a result in $\sL^p$, then we will show that if the
problem is formulated in $\sD$ it is possible to approximate it with a
sequence of $\sL^p$ problems; the solutions to such approximating
problems will be finally shown to converge to a function that solves
the Komogorov backward PDE in the sense of Definition \ref{defsolD}.

All the above discussion about the meaning of Kolmogorov equation
applies verbatim to the space $\sL^p$. A solution in $\sL^p$ is
defined in a straightforward way as follows.
\begin{definition}
\label{defsolL}
Given $\Phi\in C_{b}^{2,\alpha} (\sL^p,\mathbb{R} ) $,
we say that
$u: [ 0,T ] \times\sL^p\rightarrow\mathbb{R}$ is a
solution of the
Kolmogorov equation in $\sL^p$ with terminal condition $\Phi$ if
\[
u\in L^{\infty} \bigl( 0,T;C_{b}^{2,\alpha} \bigl(
\sL^p,\mathbb{R}%
\bigr) \bigr) \cap C \bigl( [ 0,T ] \times
\sL^p,\mathbb{R}%
\bigr)
\]
$u (\cdot,y )$ is Lipschitz for any $y\in\Dom
(A )$ and satisfies identity (\ref{eqPKolmogorov}) for every
$t\in[ 0,T ] $ and $y\in\Dom(A ) $, with the
duality terms understood with respect to the topology of $\sL^p$.
\end{definition}
%
\section{Solution in $\sL^p$}
\label{secLp}
The choice to work in a general $\sL^p$ space instead of working with
the Hilbert space $\sL^2$ could seem unjustified at first sight. As
long as solving Kolmogorov equation in $\sL^p$ is only a step toward
solving it in $\sD$ through approximations it would be enough to
develop the theory in $\sL^2$, where the results needed are well
known. Nevertheless, we give and prove here this more general statement
for $\sL^p$ spaces for some reasons. First, the proof shows a method
to obtain this kind of result without actually using a It\^o-type
formula, but only a Taylor expansion; the difference is tiny but it
allows to work in spaces where there is no It\^o formula to apply.
Second, the proof points out where a direct argument of this kind
(which is essentially the classical scheme for these results) fails.
Last, also the easiest examples do not behave well in $\sL^2$ but they
can be regular enough in some $\sL^p$ instead (see Examples \ref
{ex1long} and \ref{ex2long} hereinafter).
Therefore, proving the result in $\sL^p$ is already enough to deal
with some examples, without the need to go further in the development
of the theory.

If $B$ satisfies Assumption \ref{assB} with $E=\sL^p$, Theorems \ref
{thmSDE}, \ref{4535454754545} and \ref{thmmarkov} yield that the SDE
%
{\renewcommand{\theequation}{14$'$ bis}
\begin{eqnarray}
\ud Y(s)= \bigl[AY(s)+B \bigl(s,Y(s) \bigr) \bigr]\udt
s+\Sigma\udt
\beta(s),
\nonumber\\[-8pt]\\[-8pt]
\eqntext{s\in[t,T],  Y(t)=y} 
\end{eqnarray}\setcounter{equation}{26}}%
admits a unique mild solution $Y^{t_0,y}(t)$ in $\sL^p$ which is
continuous in time, $C^{2,\alpha}_b$ with respect to $y$ and has the
Markov property.
\begin{theorem}
\label{thmLp}
Let $\Phi:\sL^p\to\bR$ be in $C^{2,\alpha}_b$ and let Assumption
\ref{assB} hold in~$\sL^p$. Then the function
\[
u ( t,y ):=\mathbb{E} \bigl[ \Phi\bigl( Y^{t,y} ( T ) \bigr) \bigr],
\qquad( t,y ) \in[ 0,T ] \times\sL^p,
\]
is a solution of the Kolmogorov equation in $\sL^p$ with terminal
condition $\Phi$.
\end{theorem}
\begin{pf}
Throughout this proof $\llVert \cdot\rrVert $ will denote the norm in
$\sL^p$
and $\langle\cdot,\cdot\rangle$ will denote duality between $\sL
^p$ and $\sL^{p^\prime}$, where $\frac{1}{p}+\frac{1}{p^\prime
}=1$.

The function $u$ has the regularity properties required by the
definition of solution: boundedness in time is straightforward, while
the fact that $\Phi$ belongs to $C_b^{2,\alpha} (\sL^p;\bR
^d )$ and the regularity properties of $Y$ with respect to the
initial data stated in Theorem \ref{4535454754545} imply, by
composition and the dominated convergence theorem, that $u$ is
continuous on $[0,T]\times\sL^p$ and $u(t,\cdot)$ is in $C^{2,\alpha
}_b (\sL^p;\bR^d )$ for every $t\in[0,T]$; the Lipschitz
property in time is a consequence of being a solution of an integral
equation where all the terms are bounded. We have thus to show that it
satisfies equation (\ref{eqPKolmogorov}). Recall that we choose $y$ in
the domain of $A$.

\begin{longlist}
\item[\textit{Step} 1.] Fix $t_0\in[0,T]$. From Markov property, for any
$t_1>t_{0}$ in $ [0,T ] $, we have
\[
u ( t_{0},y ) =\mathbb{E} \bigl[ u \bigl( t_1,Y^{t_{0},y}
( t_1 ) \bigr) \bigr]
\]
because
\begin{eqnarray*}
\mathbb{E} \bigl[ \Phi\bigl( Y^{t_{0},y} ( T ) \bigr) \bigr] & =&\mathbb{E}
\bigl[ \mathbb{E} \bigl[ \Phi\bigl( Y^{t_{0},y} ( T ) \bigr) \mid
Y^{t_{0},y} ( t_1 ) \bigr] \bigr]
\\
& =&\mathbb{E} \bigl[ \mathbb{E} \bigl[ \Phi\bigl( Y^{t_1,w} (T ) \bigr)
\bigr] _{w=Y^{t_{0},y} ( t_1 ) } \bigr]=\mathbb{E} \bigl[ u \bigl(
t_1,Y^{t_{0},y}
( t_1 ) \bigr) \bigr].
\end{eqnarray*}
From Taylor formula applied to the function $y\mapsto u (
t,y )$, we have
\begin{eqnarray*}
&& u \bigl( t_1,Y^{t_{0},y} ( t_1 ) \bigr) -u \bigl(
t_1,e^{ (t_1-t_{0} ) A}y \bigr)
\\
&&\qquad = \bigl\langle Du \bigl( t_1,e^{ (t_1-t_{0} ) A}y \bigr)
,Y^{t_{0},y} ( t_1 ) -e^{ (t_1-t_{0} ) A}y \bigr\rangle
\\
&&\quad\qquad{}  +\tfrac{1}{2}D^{2}u \bigl( t_1,e^{ ( t_1-t_{0}
) A}y
\bigr) \bigl(Y^{t_{0},y} ( t_1 ) -e^{ (
t_1-t_{0} ) A}y,Y^{t_{0},y}
( t_1 ) -e^{ (
t_1-t_{0} ) A}y \bigr)
\\
&&\quad\qquad{} +\tfrac{1}{2}r_{u ( t_1,\cdot) } \bigl( Y^{t_{0},y}
(t_1 ), e^{ ( t_1-t_{0} ) A}y \bigr),
\end{eqnarray*}
where
\begin{eqnarray*}
&& \bigl\llvert r_{u ( t_1,\cdot) } \bigl( Y^{t_{0},y} ( t_1
),e^{ ( t_1-t_{0} ) A}y \bigr) \bigr\rrvert
\\
&&\qquad \leq\mathfrak{w} \bigl(D^{2}u ( t_1,\cdot), \bigl\llVert
Y^{t_{0},y} ( t_1 )-e^{ ( t_1-t_{0} ) A}y\bigr\rrVert\bigr) \bigl
\llVert Y^{t_{0},y} ( t_1 ) -e^{ (
t_1-t_{0} ) A}y\bigr\rrVert
^{2}
\end{eqnarray*}
(for the definitions of $r$ and $\mathfrak{w}$ see Remark \ref{RemTaylor}).
Due to the $C_{b}^{2,\alpha} ( \sL^p,\mathbb{R} ) $-property,
uniform in time, we have
\[
\bigl\llvert r_{u ( t_1,\cdot) } \bigl( Y^{t_{0},y} ( t_1
),e^{ ( t_1-t_{0} ) A}y \bigr) \bigr\rrvert\leq M\bigl\llVert
Y^{t_{0},y} (
t_1 ) -e^{ ( t_1-t_{0}
) A}y\bigr\rrVert^{2+\alpha}.
\]
Recall that
\begin{eqnarray*}
Y^{t_{0},y} ( t_1 ) -e^{ ( t_1-t_{0} ) A}y & =&F^{t_{0}}
\bigl( t_1,Y^{t_0,y} \bigr) +Z^{t_{0}} ( t_1
),
\\
F^{t_{0}} \bigl( t_1,Y^{t_0,y} \bigr) & =&\int
_{t_{0}}^{t_1}e^{ (
t_1-s )A}B \bigl(
s,Y^{t_{0},y} ( s ) \bigr)\udt s
\end{eqnarray*}
and
\begin{eqnarray*}
\mathbb{E} \bigl[ Z^{t_{0}} ( t_1 ) \bigr] & =&0,
\\
\mathbb{E} \bigl[ \bigl\llVert Z^{t_{0}} ( t_1 ) \bigr\rrVert
^{4} \bigr] & \leq& C_{Z}^{4} (
t_1-t_{0} )^{2},
\\
\bigl\llVert F^{t_{0}} \bigl( t_1,Y^{t_0,y} \bigr)
\bigr\rrVert& \leq& C{\llVert B\rrVert}_\infty( t_1-t_{0}
),
\end{eqnarray*}
where ${\llVert B\rrVert }_\infty=\sup_t\sup_y\llVert
B(t,y)\rrVert $.

Hence, recalling $u ( t_{0},y ) =\mathbb{E} [ u
(t_1,Y^{t_{0},y} ( t_1 ) ) ] $,
\begin{eqnarray*}
&& u ( t_{0},y ) -u \bigl( t_1,e^{ ( t_1-t_{0} )
A}y \bigr)
\\
&&\qquad  = \bigl\langle Du \bigl( t_1,e^{ ( t_1-t_{0} ) A}y \bigr),\mathbb{E}
\bigl[ F^{t_{0}} \bigl( t_1,Y^{t_0,y} \bigr) \bigr]
\bigr\rangle
\\
&&\quad\qquad{}  +\tfrac{1}{2}\mathbb{E} \bigl[ D^{2}u \bigl(
t_1,e^{
( t_1-t_{0} )A}y \bigr) \bigl( F^{t_{0}} \bigl(
t_1,Y^{t_0,y} \bigr) +Z^{t_{0}} ( t_1 ),
\\
&&\hspace*{2pt}\qquad\quad
F^{t_0} \bigl( t_1,Y^{t_0,y} \bigr)
+Z^{t_{0}} ( t_1 ) \bigr) \bigr]
\\
&&\quad\qquad{} +\tfrac{1}{2}\mathbb{E} \bigl[ r_{u ( t_1,\cdot
) }
\bigl(Y^{t_{0},y} ( t_1 ), e^{ (
t_1-t_{0} ) A}y \bigr) \bigr].
\end{eqnarray*}

\item[\textit{Step} 2.] Now let us explain the strategy. Given $t\in
[0,T ] $, taken a sequence of partitions $\pi_{n}$ of $
[t,T ] $, of the form $t=t_{1}^{n}\leq\cdots\leq t_{k_{n}+1}^{n}=T$
of $ [ t,T ] $ with  $\llvert \pi_{n}\rrvert
\rightarrow0$, we take $t_{0}=t_{i}^{n}$ and $t_1=t_{i+1}^{n}$ in the
previous identity and sum over the partition $\pi_{n}$ to get
\[
u ( t,y ) -\Phi( y ) +I_{n}^{ ( 1
)}=I_{n}^{ ( 2 ) }+I_{n}^{ ( 3 ) }+I_{n}^{
(4 ) },
\]
where
\begin{eqnarray*}
I_{n}^{ ( 1 ) }&:=&\sum_{i=1}^{k_{n}}
\bigl( u \bigl( t_{i+1}^{n},y \bigr) -u \bigl(
t_{i+1}^{n},e^{ (
t_{i+1}^{n}-t_{i}^{n} )A}y \bigr) \bigr),
\\
I_{n}^{ ( 2 ) }&:=&\sum_{i=1}^{k_{n}}
\bigl\langle Du \bigl(t_{i+1}^{n},e^{ ( t_{i+1}^{n}-t_{i}^{n} ) A}y
\bigr),
\mathbb{E} \bigl[ F^{t_{i}^{n}} \bigl( t_{i+1}^n,Y^{t_{i}^n,y}
\bigr) \bigr] \bigr\rangle,
\\
I_{n}^{ ( 3 ) }&:=&\frac{1}{2}\sum
_{i=1}^{k_{n}}\mathbb{E} \bigl[D^{2}u \bigl(
t_{i+1}^{n},e^{ (
t_{i+1}^{n}-t_{i}^{n} ) A}y \bigr)
\\
&&{}\times \bigl( F^{t_{i}^{n}} \bigl( t_{i+1}^{n},Y^{t_i^n,y}
\bigr)+Z^{t_{i}^{n}} \bigl(t_{i+1}^{n} \bigr)
,F^{t_{i}^{n}} \bigl( t_{i+1}^{n},Y^{t_i^n,y} \bigr)
+Z^{t_{i}^{n}} \bigl( t_{i+1}^{n} \bigr) \bigr) \bigr],
\\
I_{n}^{ ( 4 ) }&:=&\frac{1}{2}\sum
_{i=1}^{k_{n}}\mathbb{E} \bigl[r_{u ( t_{i+1}^{n},\cdot) } \bigl(
Y^{t_{i}^{n},y} \bigl(t_{i+1}^{n} \bigr), e^{ (
t_{i+1}^{n}-t_{i}^{n} ) A}y
\bigr) \bigr].
\end{eqnarray*}
We want to show that:
\begin{longlist}[(III)]
\item[(I)] $\lim_{n\rightarrow\infty}I_{n}^{ (
1 ) }=-\int_{t}^{T} \langle Du (
s,y ), Ay \rangle\udt s\mbox{ if }y\in\Dom(
A )$,

\item[(II)] $\lim_{n\rightarrow\infty}I_{n}^{ (
2 ) }=\int_{t}^{T} \langle Du (
s,y ), B ( s,y ) \rangle\udt s$,

\item[(III)] $\lim_{n\rightarrow\infty}I_{n}^{ (
3 ) }=\frac{1}{2}\int_{t}^{T}\sum
_{j=1}^d\sigma_j^{2}D^{2}u ( s,y ) ( e_j,e_j )\udt s$,

\item[(IV)] $\lim_{n\rightarrow\infty}I_{n}^{ (
4 ) }=0$.
\end{longlist}

\item[\textit{Step} 3.] We have, for $y\in\Dom( A ) $ (in this
case $\frac{\ud}{\ud t}e^{tA}y=Ae^{tA}y$)
\begin{eqnarray*}
&& \sum_i^{k_n}u \bigl(
t_{i+1}^{n}, y \bigr) -u \bigl( t_{i+1}^{n},e^{ (
t_{i+1}^{n}-t_{i}^{n} ) A}y
\bigr)
\\
&&\qquad =-\sum_i^{k_n}\int
_{0}^{t_{i+1}^{n}-t_{i}^{n}} \bigl\langle Du \bigl(
t_{i+1}^{n},e^{sA}y \bigr), Ae^{sA}y
\bigr\rangle\udt s
\\
&&\qquad  =-\sum_i^{k_n}\int
_{t_{i}^{n}}^{t_{i+1}^{n}} \bigl\langle Du \bigl(
t_{i+1}^{n},e^{
( s-t_{i}^{n} ) A}y \bigr), Ae^{ ( s-t_{i}^{n} )
A}y
\bigr\rangle\udt s
\\
&&\qquad  =-\int_t^T\sum_i^{k_n}
\bigl\langle Du \bigl(t_{i+1}^n,e^{(s-t^n_i)A}y
\bigr),Ae^{(s-t_i^n)A}y \bigr\rangle\ind_{[t_i^n,t_{i+1}^n]}(s)\udt s.
\end{eqnarray*}
The semigroup $e^{tA}$ is strongly continuous in $\sL^p$ therefore it
converges to the identity as $t$ goes to $0$; hence, since $y$ is
fixed, taking the limit in $n$ yields (I) applying the dominated
convergence theorem.

\item[\textit{Step} 4.]
By standard properties of the Bochner integral, we have
\begin{eqnarray*}
&& \sum_{i=1}^{k_n} \biggl\langle Du
\bigl(t_{i+1}^n,e^{(t^n_{i+1}-t^n_i)A}y \bigr),\bE\int
_{t_i^n}^{t^n_{i+1}}e^{(t^n_{i+1}-s)A}B \bigl(s,Y^{t_i^n,y}(s)
\bigr)\udt s \biggr\rangle
\\
&&\qquad =\sum_{i=1}^{k_n}\bE\int
_{t^n_i}^{t^n_{i+1}} \bigl\langle Du \bigl(t^n_{i+1},e^{(t^n_{i+1}-t^n_i)A}y
\bigr),e^{(t^n_{i+1}-s)A}B \bigl(s,Y^{t_i^n,y}(s) \bigr) \bigr\rangle\udt s
\\
&&\qquad =\bE\int_t^T\sum
_{i=1}^{k_n} \bigl\langle Du \bigl(t^n_{i+1},e^{(t^n_{i+1}-t^n_i)A}y
\bigr),e^{(t^n_{i+1}-s)A}B \bigl(s,Y^{t_i^n,y}(s) \bigr) \bigr\rangle
\ind_{[t^n_i,t^n_{i+1}]}(s)\udt s;
\end{eqnarray*}
now arguing as in the previous step it's easy to prove that this
quantity converges to
\[
\int_t^T \bigl\langle Du(s,y),B(s,y) \bigr
\rangle\udt s.
\]
\item[\textit{Step} 5.]
First, split each of the addends appearing in $I_n^{(3)}$ as follows:
\begin{eqnarray*}
&& D^{2}u  \bigl( t_{i+1}^{n},e^{ ( t_{i+1}^{n}-t_{i}^{n} )
A}y
\bigr)
\\
&&\quad{}\times \bigl( F^{t_{i}^{n}} \bigl( t_{i+1}^{n},Y^{t_i^n,y}
\bigr) +Z^{t_{i}^{n}} \bigl(t_{i+1}^{n} \bigr),
F^{t_{i}^{n}} \bigl( t_{i+1}^{n},Y^{t_i^n,y}
\bigr) +Z^{t_{i}^{n}} \bigl( t_{i+1}^{n} \bigr) \bigr)
\\
&&\qquad =D^{2}u \bigl( t_{i+1}^{n},e^{ ( t_{i+1}^{n}-t_{i}^{n} )
A}y
\bigr) \bigl( F^{t_{i}^{n}} \bigl( t_{i+1}^{n},Y^{t_i^n,y}
\bigr),F^{t_{i}^{n}} \bigl( t_{i+1}^{n},Y^{t_i^n,y}
\bigr) \bigr)
\\
&&\quad\qquad{} +D^{2}u \bigl( t_{i+1}^{n},e^{ (
t_{i+1}^{n}-t_{i}^{n} ) A}y
\bigr) \bigl( F^{t_{i}^{n}} \bigl( t_{i+1}^{n},Y^{t_i^n,y}
\bigr),Z^{t_{i}^{n}} \bigl( t_{i+1}^{n} \bigr) \bigr)
\\
&&\quad\qquad{}+D^{2}u \bigl( t_{i+1}^{n},e^{ (
t_{i+1}^{n}-t_{i}^{n} ) A}y
\bigr) \bigl(Z^{t_{i}^{n}} \bigl(t_{i+1}^{n} \bigr)
,F^{t_{i}^{n}} \bigl( t_{i+1}^{n},Y^{t_i^n,y} \bigr)
\bigr)
\\
&&\quad\qquad{}+D^{2}u \bigl( t_{i+1}^{n},e^{ (
t_{i+1}^{n}-t_{i}^{n} ) A}y
\bigr) \bigl(Z^{t_{i}^{n}} \bigl(t_{i+1}^{n} \bigr)
,Z^{t_{i}^{n}} \bigl( t_{i+1}^{n} \bigr) \bigr).
\end{eqnarray*}
Let us give the main estimates. We have
\begin{eqnarray*}
&& \bigl\llvert\mathbb{E} \bigl[ D^{2}u \bigl( t,e^{ ( t-t_{0} )
A}y
\bigr) \bigl( F^{t_{0}} \bigl( t,Y^{t_0,y} \bigr), F^{t_{0}}
\bigl( t,Y^{t_0,y} \bigr) \bigr) \bigr] \bigr\rrvert
\\
&&\qquad  \leq {\bigl\llVert D^{2}u\bigr\rrVert}_\infty\mathbb{E}
\bigl[ \bigl\llVert F^{t_{0}} \bigl( t,Y^{t_0,y} \bigr) \bigr\rrVert
^{2} \bigr]
\\
&&\qquad  \leq \bigl\llVert D^{2}u\bigr\rrVert_{\infty}C^{2}{
\llVert B\rrVert}_{\infty}^{2} ( t-t_{0} )
^{2}
\end{eqnarray*}
and
\begin{eqnarray*}
&& \bigl\llvert\mathbb{E} \bigl[ D^{2}u \bigl( t,e^{ ( t-t_{0} )
A}y
\bigr) \bigl( F^{t_{0}} \bigl( t,Y^{t_0,y} \bigr), Z^{t_{0}}
( t ) \bigr) \bigr] \bigr\rrvert
\\
&&\qquad  \leq \bigl\llVert D^{2}u\bigr\rrVert_\infty\mathbb{E}
\bigl[ \bigl\llVert F^{t_{0}} \bigl( t,Y^{t_0,y} \bigr) \bigr\rrVert
^{2} \bigr] ^{1/2}\mathbb{E} \bigl[ \bigl\llVert
Z^{t_{0}} ( t )\bigr\rrVert^{2} \bigr] ^{1/2}
\\
&&\qquad  \leq \bigl\llVert D^{2}u\bigr\rrVert_\infty C\cdot
C_{Z}{\llVert B\rrVert}_\infty( t-t_{0} )
^{3/2},
\end{eqnarray*}
where we have set
\[
{\bigl\llVert D^2u\bigr\rrVert}_\infty=\sup
_t\sup_y{\bigl\llVert
D^2u(t,y)\bigr\rrVert}_{L (E,E;E )},
\]
hence the first three terms give no contribution when summing up over
$i$, because they are estimated by a power of $t_{i+1}-t_i$ greater
than $1$.

Therefore, it remains to show that the term
%
\begin{equation}
\label{eqIn^3bis} \sum_{i=1}^{k_n}\bE
\bigl[D^{2}u \bigl( t_{i+1}^{n},e^{ (
t_{i+1}^{n}-t_{i}^{n} ) A}y
\bigr) \bigl(Z^{t_{i}^{n}} \bigl(t_{i+1}^{n} \bigr)
,Z^{t_{i}^{n}} \bigl( t_{i+1}^{n} \bigr) \bigr) \bigr]
\end{equation}
converges to
\[
\int_{t_0}^t\sigma^2
D^2u(s,y) (e,e)\udt s.
\]
To this aim, we recall that
\begin{eqnarray*}
Z^{t_i^n} \bigl(t_{i+1}^n \bigr)&=&\int
_{t_i^n}^{t_{i+1}^n}e^{
(t_{i+1}^n-r )A}\pmatrix{\sigma\udt  W(r)
\cr
0}
\\
&=&\pmatrix{ \sigma\bigl(W \bigl(t_{i+1}^n \bigr)-W
\bigl(t_i^n \bigr) \bigr)
\vspace*{3pt}\cr
\sigma\bigl(W \bigl(
\bigl(t_{i+1}^n+\cdot\bigr)\vee t_i^n
\bigr)-W \bigl(t_i^n \bigr) \bigr)}
\\
&=:& \pmatrix{Z^i_0
\vspace*{3pt}\cr
Z^i_1}.
\end{eqnarray*}
We split again (\ref{eqIn^3bis}) into
\begin{eqnarray*}
&& \sum_{i=1}^{k_n}\bE\biggl[D^{2}u
\bigl( t_{i+1}^{n},e^{ (
t_{i+1}^{n}-t_{i}^{n} ) A}y \bigr) \biggl(
\pmatrix{Z^i_0
\cr
0},\pmatrix{Z^i_0
\cr
0} \biggr)
\\
&&\qquad{} +D^{2}u \bigl( t_{i+1}^{n},e^{ (
t_{i+1}^{n}-t_{i}^{n} ) A}y
\bigr) \biggl(\pmatrix{Z^i_0
\cr
0},\pmatrix{0
\cr
Z^i_1} \biggr)
\\
&&\qquad{} +D^2u \bigl( t_{i+1}^{n},e^{ (
t_{i+1}^{n}-t_{i}^{n} ) A}y
\bigr) \biggl(\pmatrix{0
\cr
Z^i_1},\pmatrix{Z^i_0
\cr
0} \biggr)
\\
&&\qquad{}+D^{2}u \bigl( t_{i+1}^{n},e^{ (
t_{i+1}^{n}-t_{i}^{n} ) A}y
\bigr) \biggl(\pmatrix{0
\cr
Z^i_1},\pmatrix{0
\cr
Z^i_1} \biggr) \biggr].
\end{eqnarray*}
For the first term we have, using It\^o isometry in $\bR^d$, that
\begin{eqnarray*}
&& \sum_{i=1}^{k_n}\bE\biggl[D^2u
\bigl(t_{i+1}^n,e^{
(t_{i+1}^n-t_i^n )A}y \bigr) \biggl(
\pmatrix{Z_0^i
\cr
0},\pmatrix{Z_0^i
\cr
0} \biggr) \biggr]
\\
&&\qquad =
\sum_{j=1}^d\sigma_j^2
\sum_{i=1}^{k_n}D^2u
\bigl(t_{i+1}^n,e^{ (t_{i+1}^n-t_i^n
)A}y \bigr)
(e_j,e_j ) \bigl(t_{i+1}^n-t_i^n
\bigr)
\end{eqnarray*}
and the right-hand side in this equation converges to
\[
\sum_{j=1}^d\sigma_j^2
\int_{t_0}^tD^2u(s,y)
(e_j,e_j)\udt s
\]
thanks to the strong continuity of $e^{tA}$.

For the second term, we can write (here $\llVert \sigma\rrVert =\max
_j\llvert \sigma_j\rrvert $)
%
\begin{eqnarray}\label{termD1}
&& \bE\biggl\llvert D^{2}u \bigl( t_{i+1}^{n},e^{ (
t_{i+1}^{n}-t_{i}^{n} ) A}y
\bigr) \biggl(\pmatrix{Z^i_0
\cr
0},\pmatrix{0
\cr
Z^i_1} \biggr)\biggr\rrvert
\\
\nonumber
&&\qquad \leq \llVert\sigma\rrVert{\bigl\llVert D^2u\bigr\rrVert
}_{\infty}
\bE\bigl[\bigl\llvert W \bigl(t^n_{i+1}
\bigr)-W \bigl(t^n_i \bigr)\bigr\rrvert\bigl\llVert W
\bigl( \bigl(t^n_{i+1}+\cdot\bigr)\vee t^n_i
\bigr)-W \bigl(t^n_i \bigr)\bigr\rrVert_{L^p}
\bigr]\hspace*{-20pt}
\\
\nonumber
&&\qquad \leq\llVert\sigma\rrVert{\bigl\llVert D^2u\bigr\rrVert
}_{\infty}\bE\biggl[\bigl\llvert W \bigl(t^n_{i+1}
\bigr)-W \bigl(t^n_i \bigr)\bigr\rrvert\biggl(\int
_0^{t^n_{i+1}-t^n_i}\bigl\llvert W(r)\bigr\rrvert
^p\udt r \biggr)^{\sfrac{1}{p}} \biggr]
\\
\nonumber
&&\qquad \leq\llVert\sigma\rrVert{\bigl\llVert D^2u\bigr\rrVert
}_{\infty} \bigl(\bE\bigl\llvert W \bigl(t^n_{i+1}
\bigr)-W \bigl(t^n_i \bigr)\bigr\rrvert^2
\bigr)^{\sfrac{1}{2}}
\\
\nonumber
&&\quad\qquad{}\times \biggl(\bE\biggl[ \biggl(\int
_0^{t^n_{i+1}-t^n_i}\bigl\llvert W(r)\bigr\rrvert
^p\udt r \biggr)^{\sfrac{2}{p}} \biggr] \biggr)^{\sfrac{1}{2}}
\\
&&\qquad \leq\llVert\sigma\rrVert{\bigl\llVert D^2u\bigr\rrVert
}_{\infty} \bigl(t^n_{i+1}-t^n_i
\bigr)^{\sfrac{1}{2}} \bigl(t^n_{i+1}-t^n_i
\bigr)^{\sfrac{1}{p}}\nonumber
\\
&&\quad\qquad{}\times \Bigl(\bE\Bigl[ \Bigl(\sup_{[0,t^n_{i+1}-t^n_i]}
\bigl(\bigl\llvert W(r)\bigr\rrvert^p \bigr) \Bigr)^{\sfrac{2}{p}}
\Bigr] \Bigr)^{\sfrac{1}{2}}\nonumber
\\
\label{termD5}
&&\qquad \leq \llVert\sigma\rrVert{\bigl\llVert D^2u\bigr
\rrVert}_{\infty} \bigl(t^n_{i+1}-t^n_i
\bigr)^{1+\sfrac{1}{p}},
\end{eqnarray}
using It\^o isometry and Burkholder--Davis--Gundy inequality, thus it
converges to zero when summing over $i$ and letting $n$ go to $\infty
$.

The third term can be shown to go to zero in the exact same way and by
the same estimates as above, one obtains that
\begin{eqnarray*}
&&\bE\biggl\llvert D^{2}u \bigl( t_{i+1}^{n},e^{ (
t_{i+1}^{n}-t_{i}^{n} ) A}y
\bigr) \biggl(\pmatrix{0
\cr
Z^i_1},\pmatrix{0
\cr
Z^i_1} \biggr)\biggr\rrvert\leq\bigl(t^n_{i+1}-t^n_i
\bigr)^{1+\sfrac{2}{p}},
\end{eqnarray*}
hence it follows that also this term gives no contribution when passing
to the limit.

\item[\textit{Step} 6.]
Since
\[
\bigl\llvert r_{u ( t,\cdot) } \bigl( Y^{t_{0},y} ( t ),e^{ ( t-t_{0}
) A}y
\bigr) \bigr\rrvert\leq M\bigl\llVert Y^{t_{0},y} ( t ) -e^{ ( t-t_{0} )
A}y
\bigr\rrVert^{2+\alpha}
\]
we have that
\begin{eqnarray*}
&& \bigl\llvert\mathbb{E} \bigl[ r_{u ( t,\cdot) } \bigl( Y^{t_{0},y} (
t )
,e^{ ( t-t_{0} ) A_{E}}y \bigr) \bigr]\bigr\rrvert
\\
&&\qquad \leq M\mathbb{E}
\bigl[ \bigl
\llVert Y^{t_{0},y} ( t ) -e^{ ( t-t_{0} ) A}y\bigr\rrVert^{2+\alpha
}
\bigr]
\\
&&\qquad  \leq K \bigl( \mathbb{E} \bigl[ \bigl\llVert F^{t_{0}} \bigl(
t,Y^{t_0,y} \bigr)\bigr\rrVert^{4} \bigr] ^{\vfrac{2+\alpha
}{4}}+
\mathbb{E} \bigl[\bigl\llVert Z^{t_{0}} ( t ) \bigr\rrVert^{4}
\bigr]^{\vfrac{2+\alpha}{4}} \bigr)
\\
&&\qquad  \leq\widetilde K ( t-t_{0} ) ^{1+\sfrac{\alpha}{2}}
\end{eqnarray*}
and from this one proves that $\lim_{n\rightarrow\infty}I_{n}^{
(4 ) }=0$.\quad\qed
\end{longlist}\noqed
\end{pf}

\begin{remark}
The point in which the above argument fails when working directly in
$\sD$ is item (III) of step 2. Indeed step 5, which is the proof
of the convergence in (III), cannot be carried out when working
with the $\sup$-norm: if we start again from (\ref{termD1}) using the
norm of $\sD$ we would end up with the estimate
%
{\renewcommand{\theequation}{29$^\prime$}
\begin{equation}
\label{eqDwrong}
\quad \bE\biggl\llvert D^{2}u \biggl(
t_{i+1}^{n},e^{ (
t_{i+1}^{n}-t_{i}^{n} ) A}y \biggr) \biggl( \pmatrix{Z^i_0
\cr
0},\pmatrix{0
\cr
Z^i_1} \biggr)\biggr\rrvert\leq{\bigl\llVert
D^2u\bigr\rrVert}_{\infty
} \bigl(t^n_{i+1}-t^n_i
\bigr)
\end{equation}\setcounter{equation}{29}}%
which is not enough to obtain the convergence to $0$ that we need.
\end{remark}
%
\section{Solution in $\mathop{\sC}\limits^{\curvearrowleft}$}
\label{secC}
We now show how to use $\sL^p$ approximations in order to obtain
classical solutions of Kolmogorov equations in the sense of Definition
\ref{defsolD}. As before, we will assume that $B$ satisfied
Assumption \ref{assB} for $E=\sD$, that is,
\[
B\in L^{\infty} \bigl( 0,T;C_{b}^{2,\alpha} (\sD,\sD)
\bigr)
\]
for some $\alpha\in( 0,1 ) $.
Suppose we have a sequence $\{J_n\}$ of linear continuous operators
from $L^p(-T,0;\bR^d)$ into $C([-T,0];\bR^d)$ such that $J_n\phi
\stackrel{n\to\infty}{\longrightarrow}\phi$ uniformly for any
$\phi
\in C([-T,0];\bR^d)$. By the Banach--Steinhaus theorem, we have that
$\sup_n\llVert J_n\rrVert _{L(C([-T,0];\bR^d);C([-T,0];\bR^d))}<\infty$;
however,\vspace*{1pt} we need a slightly stronger property, namely that $\llVert
J_nf\rrVert
_\infty\leq C_J\llVert f\rrVert _\infty$ for all $f$ with at most one jump,
uniformly in $n$. Then we can define the sequence of operators
%
\begin{eqnarray}
\label{eqBn}
\nonumber
& \displaystyle B_n\colon[0,T]\times\sL^p\to
\bR^d\times\{0\},&
\nonumber\\[-8pt]\\[-8pt]\nonumber
& \displaystyle B_n(t,y)=B_n \biggl(t,\pmatrix{x
\cr
\phi} \biggr)=B_n(t,x,
\phi):=B (t,x,J_n\phi).&
\end{eqnarray}
We will often write $J_n {x \choose \phi}$ for ${x \choose
J_n\phi}$ in
the sequel.\vspace*{2pt}

It can be easily proved that if $B$ satisfies Assumption \ref{assB}
in $\sD$ then for every $n$ the operator $B_n$ satisfies Assumption
\ref{assB} both in $\sD$ and in $\sL^p$.
Thus, if we consider the approximated SDE
%
\begin{equation}
\label{eqDYshortn} \ud Y_n(t)=AY_n(t)\udt t+B_n
\bigl(t,Y_n(t)\bigr)\udt t+\Sigma\udt\tilde\beta(t),\qquad
Y_n(s)=y\in\sL^p
\end{equation}
by Theorem \ref{thmSDE} it admits a unique path by path mild solution
$Y_n^{s,y}$ such that, thanks to Theorem \ref{4535454754545}, the map
$t\mapsto Y^{s,y}_n(t)$ is in $C^{2,\alpha}_b$.
Suppose also we are given a terminal condition $\Phi\colon\sD\to\bR
$ for the backward Kolmogorov equation (\ref{eqPKolmogorov})
associated to the original problem with $B$; approximations $\Phi_n$
can be defined in the exact same way. We have then a sequence of
approximated backward Kolmogorov equations in $\sL^p$, namely
\begin{eqnarray}
\label{eqPKolmogorovn} u_n ( t,y ) -\Phi( y ) &=&\int_{t}^{T}
\bigl\langle Du_n ( s,y ), Ay+B_n ( s,y ) \bigr\rangle\udt
s
\nonumber\\[-8pt]\\[-8pt]\nonumber
&&{} +\frac{1}{2}\int_{t}^{T}\sum
_{j=1}^d\sigma_j^2D^{2}u_n
( s,y ) (e_j,e_j)\udt s
\end{eqnarray}
with terminal condition $u_n(T,\cdot)=\Phi_n$.
Theorem \ref{thmLp} yields in fact that for each $n$ the function
%
\begin{equation}
\label{eqPsoln} u_n(s,y)=\bE\bigl[\Phi_n
\bigl(Y_n^{s,y}(T) \bigr) \bigr]
\end{equation}
is a solution to equation (\ref{eqPKolmogorovn}) in $\sL^p$. If we
choose the initial condition $y$ in the space $\aC$ then
$Y_n^{s,y}(t)\in\aC$ as well for every $n$ and every $t\in[s,T]$.

An example of a sequence $\{J_n\}$ satisfying the required properties
can be constructed as follows: for any $\epsilon\in(0,\frac
{T}{2} )$ define a function $\tau_\epsilon:[-T,0]\to[-T,0]$ as
\[
\tau_\epsilon(x)= %
\cases{ -T+\epsilon, &\quad if $x\in[-T,-T+
\epsilon]$,
\cr
x, &\quad if $x\in[-T+\epsilon,-\epsilon]$,
\cr
-\epsilon, &\quad if
$x\in[-\epsilon,0]$.}
\]
Then choose any $C^{\infty}(\bR;\bR)$ function $\rho$ such that
$\llVert \rho\rrVert _1=1$, $0\leq\rho\leq1$ and $\supp(\rho
)\subseteq
[-1,1]$ and define a sequence $ \{\rho_n \}$ of mollifiers
by $\rho_n(x):=n\rho(nx)$. Finally set, for any $\phi\in
L^1(-T,0;\bR^d)$
%
\begin{equation}
\label{eqJn} J_n\phi(x):=\int_{-T}^0
\rho_n \bigl(\tau_{\sfrac{1}{n}}(x)-y \bigr)\phi(y)\udt y.
\end{equation}
We will need one further assumption, together with the required
properties for $J_n$ that we write again for future reference.
\begin{definition}
\label{defphi}
Let $F$ be a Banach space, $R\colon\sD\to F$ a twice Fr\'echet
differentiable function and $\Gamma\subseteq\sD$. We say that $R$
has \emph{one-jump-continuous Fr\'echet differentials of first and
second order on $\Gamma$} if there exists a sequence of linear
continuous operators $J_n:L^p(-T,0;\bR^d)\to C([-T,0];\bR^d)$ such
that $J_n\phi\stackrel{n\to\infty}{\longrightarrow}\phi$ uniformly
for any $\phi\in C([-T,0];\bR^d)$, $\sup_n\llVert J_n\phi\rrVert
_\infty\leq C_J\llVert \phi\rrVert _\infty$ for every $\phi$ that has
at most one jump and is continuous elsewhere and such that for every
$y\in\Gamma$ and for almost every $a\in[-T,0]$ the following hold:
\begin{eqnarray*}
& \displaystyle DR(y)J_n\pmatrix{1
\cr
\ind_{[a,0)}}\longrightarrow DR(y)
\pmatrix{1
\cr
\ind_{[a,0)}},&
\\
&\displaystyle D^2R(y) \biggl(J_n\pmatrix{1
\cr
\ind_{[a,0)}}-
\pmatrix{1
\cr
\ind_{[a,0)}},\pmatrix{1
\cr
\ind_{[a,0)}} \biggr)
\longrightarrow 0,&
\\
&\displaystyle D^2R(y) \biggl(\pmatrix{1
\cr
\ind_{[a,0)}},J_n
\pmatrix{1
\cr
\ind_{[a,0)}}-\pmatrix{1
\cr
\ind_{[a,0)}} \biggr)
\longrightarrow 0,&
\\
&\displaystyle D^2R(y) \biggl(J_n\pmatrix{1
\cr
\ind_{[a,0)}}-
\pmatrix{1
\cr
\ind_{[a,0)}},J_n\pmatrix{1
\cr
\ind_{[a,0)}}-\pmatrix{1
\cr
\ind_{[a,0)}} \biggr)\longrightarrow 0,&
\end{eqnarray*}
where we adopt the convention that ${1 \choose \ind
_{[a,0)}}={1 \choose 0}$ when $a=0$.
\end{definition}
We will call a sequence $ \{J_n \}$ as above a \emph
{smoothing sequence}.
\begin{assumption}
\label{assdphi}
For any $r\in[0,T]$, $B(r,\cdot)$ and $\Phi$ have
one-jump-conti\-nuous Fr\'echet differentials of first and second order
on $\aC$ and the smoothing sequence of $B$ does not depend on $r$.
\end{assumption}
\begin{remark}
Assumption \ref{assdphi} implies that the same set of properties
holds if we substitute ${1 \choose \ind_{[a,0]}}$ with any element
$q={\psi(0) \choose \psi}\in\sD_{-a}$, that is, it has at most one
jump and no other discontinuities; this happens by linearity, because
any such $\psi$ is the sum of a continuous function and an indicator function.
\end{remark}
We state and prove now the main result in this work.
\begin{theorem}
\label{thmmain}
Let $\Phi\in C^{2,\alpha}_b (\sD,\bR)$ be given and let
Assumption \ref{assB} hold for $E=\sD$. Under Assumption \ref
{assdphi}, the function $u:[0,T]\times\sD\to\bR$ given by
%
\begin{equation}
\label{eqsolu} u(t,y)=\bE\bigl[\Phi\bigl(Y^{t,y}(T) \bigr) \bigr],
\end{equation}
where $Y^{t,y}$ is the solution to equation (\ref{eqSDEHilbert})
in $\sD$, is a classical solution of the Kolmogorov equation with
terminal condition $\Phi$, that is, for every $(t,y)\in[0,T]\times
\Dom(A_{\aC} )$ it satisfies identity
{\renewcommand{\theequation}{26}
\begin{eqnarray}
u ( t,y ) -\Phi( y ) &=&\int_{t}^{T}
\bigl\langle Du ( s,y ), Ay+B ( s,y ) \bigr\rangle\udt s
\nonumber\\[-8pt]\\[-8pt]\nonumber
&&{} + \frac{1}{2}\int_{t}^{T}\sum
_{j=1}^d\sigma_j^2
D^2u(s,y) (e_j,e_j)\udt s.
\end{eqnarray}\setcounter{equation}{35}}%
\end{theorem}
\begin{pf}
Throughout this proof, $\llVert \cdot\rrVert $ will denote the norm of
$\sD$.
Sometimes we will write ${\llVert y\rrVert }_{\aC}$ to stress the
fact that
$y$ belongs to $\aC$. The duality $\langle\cdot,\cdot\rangle$ will
be always intended between $\sD^\prime$ and $\sD$.
We suppose here for simplicity that we can choose the same sequence
$ \{J_n \}$ for $B$ and $\Phi$ in Assumption \ref
{assdphi}; this does not turn in a loss of generality and the proof
can be carried on in the same way also when the two smoothing sequences
are different.

Using that smoothing sequence define $B_n$, $\Phi_n$, $Y_n$ and $u_n$
as above. The proof will be divided into some steps that will prove the
following: for $y\in\Dom(A_{\aC} )$:\vspace*{-6pt}
\begin{itemize}
\item[$\diamond$] $Y_n^{s,y}(t)\to Y^{s,y}(t)$ in $\aC$ for every
$t$ uniformly in $\omega$;
\item[$\diamond$] $u_n(s,y)\to u(s,y)=\bE[\Phi
(Y^{s,y}(T) ) ]$ for every $s$ pointwise in $y$;
\item[$\diamond$] equation (\ref{eqPKolmogorovn}) converges to
equation (\ref{eqPKolmogorov}) for any $t\in[0,T]$.
\end{itemize}

\begin{longlist}
\item[\textit{Step} 1.]
Fix $\omega\in\Omega_0$. We first need to compute
\begin{eqnarray}
\nonumber
&& \bigl\llVert Y^{s,y}_n(t)-Y^{s,y}(t)\bigr
\rrVert_{\aC}
\\
\nonumber
&&\qquad =\biggl\llVert\int_s^t
e^{(t-r)A}B_n \bigl(r,Y^{s,y}_n(r) \bigr)\udt r-\int_s^t e^{(t-r)A}B
\bigl(r,Y^{s,y}(r) \bigr)\udt r\biggr\rrVert_{\aC}
\\
\label{eqYn-Ya}
&&\qquad \leq\biggl\llVert\int_s^t
e^{(t-r)A}B_n \bigl(r,Y^{s,y}(r) \bigr)\udt r-\int
_s^t e^{(t-r)A}B \bigl(r,Y^{s,y}(r)
\bigr)\udt r\biggr\rrVert_{\aC}
\\
\label{eqYn-Yb}&&\quad\qquad{}+\biggl\llVert\int_s^t
e^{(t-r)A}B_n \bigl(r,Y^{s,y}_n(r) \bigr)\udt r-\int_s^t e^{(t-r)A}B_n
\bigl(r,Y^{s,y}(r) \bigr)\udt r\biggr\rrVert_{\aC
}.
\end{eqnarray}
For the term (\ref{eqYn-Ya}), recall that
\[
e^{(t-r)A}B_n \bigl(r,Y^{s,y}(r)
\bigr)=e^{(t-r)A}B \bigl(r,J_nY^{s,y}(r) \bigr)
\]
and that, thanks to the properties of $J_n$,
\[
J_nY^{s,y}(r)\stackrel{n\to\infty} {\longrightarrow}Y^{s,y}(r)
\]
in $\aC$, hence by continuity of $B$
%
\begin{equation}
\label{eqconvBn} B \bigl(r,J_nY^{s,y}(r) \bigr)
\longrightarrow B \bigl(r,Y^{s,y}(r) \bigr)
\end{equation}
pointwise as functions of $r$. Since $B$ is uniformly bounded in $r\in
[s,t]$, by the dominated convergence theorem
\[
\lim_{n\to\infty}\int_s^t
e^{(t-r)A}B_n \bigl(r,Y^{s,y}(r) \bigr)\udt r=\int
_s^t e^{(t-r)A}B \bigl(r,Y^{s,y}(r)
\bigr)\udt r;
\]
hence for any $\epsilon>0$
%
\begin{equation}
\label{eqprimotermine} \biggl\llVert\int_s^t
e^{(t-r)A}B_n \bigl(r,Y^{s,y}(r) \bigr)\udt r-\int
_s^t e^{(t-r)A}B \bigl(r,Y^{s,y}(r)
\bigr)\udt r\biggr\rrVert_{\aC}<\epsilon
\end{equation}
for $n$ big enough. Consider now (\ref{eqYn-Yb}):
\begin{eqnarray*}
&& \biggl\llVert\int_s^t e^{(t-r)A}B_n
\bigl(r,Y^{s,y}_n(r) \bigr)\udt r-\int_s^t
e^{(t-r)A}B_n \bigl(r,Y^{s,y}(r) \bigr)\udt r\biggr
\rrVert_{\aC}
\\
&&\qquad \leq C\int_s^t\bigl\llVert B
\bigl(r,J_nY^{s,y}_n(r) \bigr)-B
\bigl(r,J_nY^{s,y}(r) \bigr)\bigr\rrVert\udt r
\\
&&\qquad \leq C\int_s^t K_B\bigl
\llVert Y^{s,y}_n(r)-Y^{s,y}(r)\bigr\rrVert\udt r
\end{eqnarray*}
because, for any $\psi\in C ([-T,0];\bR^d )$, $\llVert
J_n\psi
\rrVert _\infty\leq C_J \llVert \psi\rrVert _\infty$ and, therefore,
$\llVert
J_ny\rrVert \leq C_J\llVert y\rrVert $.
Hence, this and (\ref{eqprimotermine}) yield, by Gronwall's lemma,
\[
\bigl\llVert Y^{s,y}_n(t)-Y^{s,y}(t)\bigr\rrVert
_{\aC}\leq\epsilon e^{TCK_B}
\]
for any $\epsilon>0$ and $n$ big enough. This implies that
$Y^{s,y}_n(t)$ converges to $Y^{s,y}(t)$ in $\aC$ for any $t$.

\item[\textit{Step} 2.]
It is now easy to deduce that $u_n(s,y)$ converges to $u(s,y)$ for any
$s$, $y\in\aC$. In fact,
\begin{eqnarray*}
&& \bigl\llvert u_n(s, y)-u(s,y)\bigr\rrvert
\\
&&\qquad \leq\bE\bigl\llvert\Phi_n \bigl(Y_n^{s,y}(T)
\bigr)-\Phi_n \bigl(Y^{s,y}(T) \bigr)\bigr\rrvert+\bE\bigl
\llvert\Phi_n \bigl(Y^{s,y}(T) \bigr)-\Phi
\bigl(Y^{s,y}(T) \bigr)\bigr\rrvert
\end{eqnarray*}
and for almost every $\omega$
\[
\bigl\llvert\Phi_n \bigl(Y^{s,y}_n(T) \bigr)-
\Phi_n \bigl(Y^{s,y}(T) \bigr)\bigr\rrvert\leq
K_\Phi\bigl\llVert Y_n^{s,y}(T)-Y^{s,y}(Y)
\bigr\rrVert
\]
and
\[
\bigl\llvert\Phi_n \bigl(Y^{s,y}(T) \bigr)-\Phi
\bigl(Y^{s,y}(T) \bigr)\bigr\rrvert\leq K_\Phi\bigl\llVert
J_nY^{s,y}(T)-Y^{s,y}(T)\bigr\rrVert,
\]
both of which are arbitrarily small for $n$ large enough; now since $B$
is bounded and we assumed that $\bE\llVert Z\rrVert ^4$ is
finite, we can apply again the dominated convergence theorem
(integrating in the variable $\omega$) to conclude this argument.

\item[\textit{Step} 3.]
We now approach the convergence of the term
\[
\bigl\langle Du_n (s,y), B_n(s,y)\bigr\rangle;
\]
it is enough to consider a generic sequence $\tilde g_n\to\tilde g$ in
$\bR^d$, to which we associate the corresponding sequence
$g_n={\tilde g_n \choose 0}\to g={\tilde g \choose 0}$ in
$\sC\subset\sD$.
From (\ref{eqPsoln}) and (\ref{eqsolu}), we have that for $h\in
\sD$
%
\begin{equation}
\label{eqDuun} \bigl\langle Du_n(s,y),h\bigr\rangle=\bE\bigl[\bigl
\langle D\Phi_n \bigl(Y_n^{s,y}(T)
\bigr),DY_n^{s,y}(T)h\bigr\rangle\bigr]
\end{equation}
and
%
\begin{equation}
\label{eqDuu} \bigl\langle Du(s,y),h\bigr\rangle=\bE\bigl[\bigl\langle D
\Phi_n \bigl(Y^{s,y}(T) \bigr),DY^{s,y}(T)h\bigr
\rangle\bigr].
\end{equation}
We remark here that the duality ${}_{\sD^{\prime}}\langle
Du_n(s,y),g_n\rangle_{\sD}$ is well defined and equals
${}_{{\sL^p}^{\prime}} \langle Du_n(s,y),g_n\rangle_{\sL^p}$; a
simple proof of this fact is the following: $u_n$ is Fr\'echet
differentiable both on $\sD$ and on $\sL^p$ and its G\^ateaux
derivatives along the directions in $\sD$ are of course the same in
$\sD$ and in $\sL^p$, therefore, also the Fr\'echet derivatives must
be equal. Now
\begin{eqnarray*}
&& \bigl\llvert\langle D u_n,g_n\rangle-\langle Du,g
\rangle\bigr\rrvert\\
&&\qquad =\bigl\llvert\langle Du_n,g_n-g
\rangle+\langle Du_n-Du,g\rangle\bigr\rrvert
\\
&&\qquad \leq\bigl\llvert\langle Du_n-Du,g\rangle\bigr\rrvert+\bigl\llvert
\langle Du_n,g_n-g\rangle\bigr\rrvert
\\
&&\qquad \leq\bE\bigl\llvert\bigl\langle D\Phi_n \bigl(Y_n^{s,y}(T)
\bigr), DY^{s,y}_n(T)g\bigr\rangle-\bigl\langle D\Phi
\bigl(Y^{s,y}(T) \bigr),DY^{s,y}(T)g\bigr\rangle\bigr\rrvert
\\
&&\quad\qquad{} +\bE\bigl\llvert\bigl\langle D\Phi_n
\bigl(Y_n^{s,y}(T) \bigr),DY_n^{s,y}(T)
(g_n-g)\bigr\rangle\bigr\rrvert
\\
&&\qquad =\bE\llvert\mA\rrvert+\bE\llvert\mathrm{B}\rrvert.
\end{eqnarray*}
We show that this last expression goes to $0$ as $n\to\infty$. We
start from $\mathrm{B}$. It is easily shown that
\[
D\Phi_n(\hat y)=D\Phi(J_n\hat y)J_n
\]
for any $\hat y\in\sD$. $D\Phi$ is bounded by assumption, whereas by
the required properties of $J_n$
\[
\bigl\llVert J_nDY_n^{s,y}(T)c\bigr\rrVert\leq
C_J \bigl\llVert DY_n^{s,y}(T)c\bigr\rrVert
\]
for any $c\in\sC$. Since the $\llVert DY_n\rrVert $'s are
uniformly bounded by a constant depending only on $e^{tA}$ and on $DB$
(see the proof of Theorem \ref{4535454754545} in the \hyperref[append]{Appendix}), we
have that the $Du_n$'s are uniformly bounded on $\sC$ as well and,
therefore, $\bE\llvert \mathrm{B}\rrvert \to0$ as $g_n\to
g$.

The term $\mA$ requires some work: from now on fix $\omega\in\Omega
_0$ and write (suppressing indexes $s$, $y$, $\omega$ and $T$)
\begin{eqnarray*}
\mA&=&\bigl\langle D\Phi_n(Y_n), DY_ng\bigr
\rangle-\bigl\langle D\Phi(Y),DYg\bigr\rangle
\\
&=&\bigl\langle D\Phi_n(Y_n),(DY_n-DY)g\bigr
\rangle+\bigl\langle D\Phi_n(Y_n)-D\Phi(Y),DYg\bigr
\rangle=\mA_1+\mA_2,
\\
\mA_2 &=&\bigl\langle D\Phi_n(Y_n)-D
\Phi_n(Y),DYg\bigr\rangle+ \bigl\langle D\Phi_n(Y)-D
\Phi(Y),DYg\bigr\rangle=\mA_{21}+\mA_{22}.
\end{eqnarray*}
Since the Lipschitz constants of $D\Phi_n$ are uniformly bounded in
$\aC$, we have that
\begin{eqnarray*}
\llvert\mA_{21}\rrvert&\leq&\bigl\llVert D\Phi_n(Y_n)-D
\Phi_n(Y)\bigr\rrVert_{\sD^{\prime}} \llVert DYg\rrVert
_{\sD
}
\\
&\leq& K_1\llVert Y_n-Y\rrVert\llVert DYg\rrVert
\end{eqnarray*}
and the last line goes to zero as $n$ goes to infinity. For $\mA
_{22}$, write
\begin{eqnarray*}
\llvert\mA_{22}\rrvert&=&\bigl\llvert\bigl\langle D\Phi
(J_nY)J_n,DYg\bigr\rangle-\bigl\langle D\Phi(Y),DYg\bigr
\rangle\bigr\rrvert
\\
&\leq&\bigl\llvert\bigl\langle D\Phi(J_nY)J_n,DYg\bigr
\rangle-\bigl\langle D\Phi(Y)J_n,DYg\bigr\rangle\bigr\rrvert
\\
&&{}+\bigl\llvert\bigl\langle D\Phi(Y)J_n,DYg\bigr
\rangle-\bigl\langle D\Phi(Y),DYg\bigr\rangle\bigr\rrvert
\\
&\leq& K_{D\Phi}\llVert J_nY-Y\rrVert\llVert DYg\rrVert+
\bigl\llvert\bigl\langle D\Phi(Y)J_n,DYg\bigr\rangle-\bigl\langle D
\Phi(Y),DYg\bigr\rangle\bigr\rrvert;
\end{eqnarray*}
the first term goes to zero by properties of $J_n$, the second one
thanks to Assumption~\ref{assdphi}: this is because from the defining
equation for $DY$ one easily sees that for any ${g \choose  0}\in
\sC$
the second component of $DYg$ has a unique discontinuity point, and our
assumption is made exactly in order to be able to control the
convergence of these terms. Now we consider $\mA_1$:
%
\begin{eqnarray}\label{eqA1}
\nonumber
&& DY_n^{s,y}(T)g-DY^{s,y}(T)g
\\
&&\qquad
= \int_s^T e^{(T-r)A}DB_n
\bigl(r,Y^{s,y}_n(r) \bigr) \bigl[DY^{s,y}_n(r)-DY^{s,y}(r)
\bigr]g\udt r
\nonumber\\[-8pt]\\[-8pt]\nonumber
&&\quad\qquad {}+\int_s^T
e^{(T-r)A} \bigl[DB_n \bigl(r,Y^{s,y}_n(r)
\bigr)-DB \bigl(r,Y^{s,y}(r) \bigr) \bigr]DY^{s,y}(r)g\udt r
\\
\nonumber
&&\qquad = \mA_{11}+\mA_{12}
\end{eqnarray}
and $\mA_{12}$ can be written as
\begin{eqnarray*}
\mA_{12}&=&\int_s^T
e^{(T-r)A} \bigl[ DB_n \bigl(r,Y^{s,y}_n(r)
\bigr)-DB_n \bigl(r,Y^{s,y}(r) \bigr) \bigr]DY^{s,y}(r)g
\udt r
\\
&&{}+\int_s^T e^{(T-r)A}
\bigl[DB_n \bigl(r,Y^{s,y}(r) \bigr)-DB \bigl(r,Y^{s,y}(r)
\bigr) \bigr] DY^{s,y}(r)g\udt r
\\
&=&\mA_{121}+\mA_{122}
\end{eqnarray*}
whence
\begin{eqnarray*}
\llVert\mA_{121}\rrVert&\leq& C\int_s^T
\bigl\llVert DY^{s,y}(r)g\bigr\rrVert\bigl\llVert DB
\bigl(r,J_nY^{s,y}_n(r) \bigr)-DB
\bigl(r,J_nY^{s,y}(r) \bigr)\bigr\rrVert\udt r
\\
&\leq& C \int_s^T \bigl\llVert
DY^{s,y}(r)g\bigr\rrVert\bigl\llVert D^2B(r,\cdot)\bigr
\rrVert\bigl\llVert J_nY^{s,y}_n(r)-
J_nY^{s,y}(r)\bigr\rrVert
\\
&\leq& C\cdot C_J \int_s^T
\bigl\llVert DY^{s,y}(r)g\bigr\rrVert\bigl\llVert D^2B(r,
\cdot)\bigr\rrVert\bigl\llVert Y^{s,y}_n(r)-Y^{s,y}(r)
\bigr\rrVert\udt r
\end{eqnarray*}
that goes to zero; for $\mA_{122}$
\begin{eqnarray*}
&& \bigl\llVert\bigl[DB_n \bigl(r,Y^{s,y}(r) \bigr)-DB
\bigl(r,Y^{s,y}(r) \bigr) \bigr]DYg\bigr\rrVert
\\
&&\qquad \leq\bigl\llVert DB \bigl(r,J_nY^{s,y}(r) \bigr)-DB
\bigl(r,Y^{s,y}(r) \bigr)\bigr\rrVert\bigl\llVert J_nDY^{s,y}(r)g
\bigr\rrVert
\\
&&\quad\qquad{}+\bigl\llVert DB \bigl(r,Y^{s,y}(r) \bigr)
\bigl[J_nDY^{s,y}(r) g-DY^{s,y}(r)g \bigr]\bigr
\rrVert
\\
&&\qquad \leq K_{DB}\bigl\llVert J_nY^{s,y}(r)-Y^{s,y}(r)
\bigr\rrVert\bigl\llVert DY^{s,y}(r)g\bigr\rrVert
\\
&&\quad\qquad{}+ \bigl\llVert DB \bigl(r,Y^{s,y}(r) \bigr)
\bigl[J_nDY^{s,y}(r) g-DY^{s,y}(r)g \bigr]\bigr
\rrVert,
\end{eqnarray*}
where the last line goes to zero thanks to Assumption \ref{assdphi}
again, and therefore $\mA_{122}$ goes to zero by the dominated
convergence theorem. From (\ref{eqA1}) and this last argument it
follows that for any fixed $\epsilon>0$
%
\begin{eqnarray}\label{eqA11}
&& \bigl\llVert DY_n(T)g-DY(T)g\bigr\rrVert
\nonumber\\[-8pt]\\[-8pt]\nonumber
&&\qquad \leq C\int
_s^T\llVert DB_n\rrVert\bigl
\llVert DY^{s,y}_n(r)g-DY^{s,y}(r)g\bigr\rrVert\udt
r+\epsilon
\end{eqnarray}
for $n$ large enough. Since $\llVert DB_n\rrVert $ is bounded
uniformly in $n$
and in $r$ we can use Gronwall's lemma to prove that $\llVert
DY_n^{s,y}(T)g-DY^{s,y}(T)g\rrVert \to0$, and since $\llVert D\Phi
_n\rrVert
$ are
uniformly bounded as well we can conclude that also $\mA_1\to0$ as
$n\to\infty$. Putting together all the pieces, we just examined we
obtain the desired convergence of $\langle Du_n,B_n\rangle$ to
$\langle Du,B\rangle$ thanks to the dominated convergence theorem (in
the variable $\omega$).

\item[\textit{Step} 4.]
All the procedures used in the previous steps apply again to treat the
convergence of the term
\[
\bigl\langle Du_n(s,y),Ay\bigr\rangle,
\]
no further passages are needed; therefore, we omit the computations and
go on to the term involving the second derivatives.

\item[\textit{Step} 5.]
We will study only the convergence of
\[
D^2u_n(s,y) (e_1,e_1)
\]
since the $\sigma_j$'s are constants and the passage from one to $d$
dimensions is trivial. We will drop the subscript $1$ in the
computations to simplify notation.
We can proceed as follows (suppressing again $s$, $y$, $\omega$ and $T$):
\begin{eqnarray*}
&& \bigl\llvert D^2u_n(s,y) (e,e)-D^2u(s,y)
(e,e)\bigr\rrvert
\\
&&\qquad \leq\bE\bigl\llvert D^2\Phi_n(Y_n)
(DY_ne,DY_ne )-D^2\Phi(Y) (DYe,DYe )\bigr
\rrvert
\\
&&\quad\qquad{} +\bE\bigl\llvert\bigl\langle D\Phi_n(Y_n),D^2Y_n(e,e)
\bigr\rangle-\bigl\langle D\Phi(Y),D^2Y(e,e)\bigr\rangle\bigr\rrvert
\\
&&\qquad =\bE\llvert\mC\rrvert+\bE\llvert\mD\rrvert.
\end{eqnarray*}
The kind of computations needed are similar to those for the terms
involving the first derivative. We first write $\mC$ (for $\omega$
fixed) as
\begin{eqnarray*}
\mC&=& \bigl[D^2\Phi_n(Y_n)
(DY_ne,DY_ne )-D^2\Phi_n(Y_n)
(DYe,DYe ) \bigr]
\\
&&{}+ \bigl[D^2\Phi_n(Y_n)
(DYe,DYe )-D^2\Phi(Y) (DYe,DYe ) \bigr]
\\
&=&\mC_1+\mC_2.
\end{eqnarray*}
For $\mC_1$, just write
\begin{eqnarray*}
\llvert\mC_1\rrvert&\leq&\bigl\llvert D^2
\Phi_n(Y_n) (DY_ne-DYe,DY_ne-DYe)
\\
&&{} +D^2\Phi_n(Y_n) (DYe,DY_ne-DYe
)+D^2\Phi_n(Y_n) (DY_ne-DYe,DYe )
\bigr\rrvert
\\
&\leq&\bigl\llVert D^2\Phi_n(Y_n)\bigr\rrVert
\bigl[\llVert DY_ne-DYe\rrVert^2+2\llVert DYe\rrVert
\llVert DY_ne-DYe\rrVert\bigr]
\end{eqnarray*}
and the last line goes to zero by the same reasoning as in $\mA_1$ and
the boundedness of $\llVert D^2\Phi_n(Y_n)\rrVert $ (uniformly in $n$).

Write $\mC_2$ as
\begin{eqnarray*}
\mC_2&=&D^2\Phi_n (Y_n ) (DYe,DYe
)-D^2\Phi(Y ) (DYe,DYe )
\\
&=&D^2\Phi(J_nY_n ) (J_nDYe,J_nDYe
)-D^2\Phi(Y ) (DYe,DYe )
\\
&=& \bigl[D^2\Phi(J_nY_n )
(J_nDYe,J_nDYe )-D^2\Phi(J_nY )
(J_nDYe,J_nDYe ) \bigr]
\\
&&{}+ \bigl[D^2\Phi(J_nY )
(J_nDYe,J_nDYe )-D^2\Phi(Y ) (DYe,DYe )
\bigr]
\\
&=&\mC_{21}+\mC_{22}.
\end{eqnarray*}
Now
\[
\mC_{21}= \bigl[D^2\Phi(J_nY_n
)-D^2\Phi(J_nY ) \bigr] (J_nDYe,J_nDYe
)
\]
hence
%
\begin{eqnarray}\label{eqestimC21}
\nonumber
\llVert\mC_{21}\rrVert&\leq&\llVert J_nDYe
\rrVert^2 \bigl\llVert D^2\Phi\bigr\rrVert
_\alpha\llVert J_nY_n-J_nY\rrVert
\nonumber\\[-8pt]\\[-8pt]\nonumber
&\leq& C_J^2\llVert DYe\rrVert\bigl
\llVert D^2\Phi\bigr\rrVert_\alpha C_J\llVert
Y_n-Y\rrVert
\end{eqnarray}
[here $\llVert D^2\Phi\rrVert _\alpha$ is the $\alpha$-H\"older norm of
$D^2\Phi$ as a map from $\sD$ into the set of bilinear forms $L
(\sD,\sD;\sD)$] which converges to zero thanks to the first
step of the proof. For $\mC_{22}$, we can write
\begin{eqnarray*}
\mC_{22}&=& \bigl[D^2\Phi(J_nY)-D^2
\Phi(Y) \bigr] (J_nDYe,J_nDYe )
\\
&&{} +D^2\Phi(Y) (J_nDYe,J_nDYe
)-D^2\Phi(Y) (DYe,DYe )
\\
&=& \bigl[D^2\Phi(J_nY)-D^2\Phi(Y) \bigr]
(J_nDYe,J_nDYe )
\\
&&{} +D^2\Phi(Y) (J_nDYe,J_nDYe-DYe
)
\\
&&{} +D^2\Phi(Y) (J_nDYe-DYe,DYe )
\\
&=& \bigl[D^2\Phi(J_nY)-D^2\Phi(Y) \bigr]
(J_nDYe,J_nDYe )
\\
&&{} +D^2\Phi(Y) (J_nDYe-DYe,J_nDYe-DYe
)
\\
&&{} +D^2\Phi(Y) (DYe,J_nDYe-DYe )
\\
&&{}+D^2\Phi(Y) (J_nDYe-DYe,J_nDYe-DYe
).
\end{eqnarray*}
Last three terms go to zero by Assumption \ref{assdphi}, while the
first one is bounded in norm by
\[
C_J\bigl\llVert D^2\Phi\bigr\rrVert_\alpha
\llVert J_nY-Y\rrVert^\alpha\llVert DYe\rrVert
^2
\]
which goes to zero since $\llVert J_nY-Y\rrVert \to0$.

We now go on with $\mD$:
\begin{eqnarray*}
\mD&=&\bigl\langle D\Phi_n(Y_n), D^2Y_n(e,e)-D^2Y(e,e)
\bigr\rangle+\bigl\langle D\Phi_n(Y_n)-D
\Phi(Y),D^2Y(e,e)\bigr\rangle
\\
&=&\mD_1+\mD_2
\end{eqnarray*}
and $\mD_2$ is easy to handle since
\[
\llvert\mD_2\rrvert\leq\bigl\llvert\bigl\langle D
\Phi_n(Y_n)-D\Phi_n(Y),D^2Y(e,e)
\bigr\rangle\bigr\rrvert+\bigl\llvert\bigl\langle D\Phi_n(Y)-D\Phi
(Y),D^2Y(e,e)\bigr\rangle\bigr\rrvert,
\]
where the first term is bounded by
\[
\bigl\llVert D^2\Phi_n\bigr\rrVert\llVert
Y_n-Y\rrVert\bigl\llVert D^2Y(e,e)\bigr\rrVert
\]
and, therefore, goes to zero as for $\mA_1$, and the second goes to
zero since $D^2Y(e,e)$ is in $\aC$ and $D\Phi_n(y)$ converge to
$D\Phi(y)$ for any $y$ as functionals on $\aC$.
Let us now rewrite the right-hand term in the bracket defining $\mD_1$ as
%
\begin{eqnarray}
\label{eqD2}
&& D^2Y_n^{s,y} (T)
(e,e)-D^2Y^{s,y}(T) (e,e)\nonumber
\\
&&\qquad =\int_s^T e^{(T-r)A}
\bigl[D^2B_n \bigl(r,Y^{s,y}_n(r)
\bigr) \bigl(DY^{s,y}_n(r) e,DY^{s,y}_n(r)e
\bigr)\nonumber
\\
&&\quad\qquad{} -D^2B \bigl(r,Y^{s,y}(r) \bigr)
\bigl(DY^{s,y}(r)e,DY^{s,y}(r)e \bigr) \bigr]\udt r
\nonumber\\[-8pt]\\[-8pt]\nonumber
&&\quad\qquad{} +\int_s^T
e^{(T-r)A} \bigl[DB_n \bigl(r,Y^{s,y}_n(r)
\bigr)D^2Y^{s,y}_n(r) (e,e)
\\
\nonumber
&&\quad\qquad{}-DB \bigl(r,Y^{s,y}(r) \bigr)D^2Y^{s,y}(r)
(e,e) \bigr]\udt r
\\
\nonumber
&&\qquad =\mD_{11}+\mD_{12}.
\end{eqnarray}
Proceeding in a way similar to before we write the integrand in $\mD
_{11}$ as a sum (suppressing also the variable $r$)
\begin{eqnarray*}
\mD_{11}&=& \bigl[D^2B_n(Y_n)
(DY_ne,DY_ne)-D^2B_n(Y_n)
(DYe,DYe) \bigr]
\\
&&{}+ \bigl[D^2B_n(Y_n)-D^2B(Y)
\bigr](DYe,DYe)
\\
&=&\mD_{111}+\mD_{112}
\end{eqnarray*}
and notice that
\begin{eqnarray*}
\llVert\mD_{111}\rrVert&=&\bigl\llVert D^2B_n(Y_n)
(DY_ne-DYe,DY_ne-DYe)
\\
&&{}+D^2B_n(Y_n)
(DY_ne-DYe,DYe) +D^2B_n(Y_n)
(DYe,DY_ne-DYe)\bigr\rrVert
\\
&\leq&\bigl\llVert D^2B_n(Y_n)\bigr\rrVert
\bigl[\llVert DY_ne-DYe\rrVert^2+2\llVert DYe\rrVert
\llVert DY_ne-DYe\rrVert\bigr]
\end{eqnarray*}
which can be treated as in $\mA_1$ since the norms $\llVert
D^2B_n (r,Y^{s,y}_n(r) )\rrVert $ are bounded uniformly in
$n$ and
$r$. $\mD_{112}$ can be treated as we did for $\mC_2$, obtaining
\begin{eqnarray*}
\mD_{112}&=& \bigl[D^2B (J_nY_n )
(J_nDYe,J_nDYe )-D^2B (J_nY )
(J_nDYe,J_nDYe ) \bigr]
\\
&&{} + \bigl[D^2B (J_nY )
(J_nDYe,J_nDYe )-D^2B (Y ) (DYe,DYe ) \bigr]
\\
&=&\mD_{1121}+\mD_{1122};
\end{eqnarray*}
an estimate analogous to (\ref{eqestimC21}) shows how to control the
term $\mD_{1121}$, while
\begin{eqnarray*}
\mD_{1122}&=& \bigl[D^2B(J_nY)-D^2B(Y)
\bigr] (J_nDYe,J_nDYe )
\\
&&{} +D^2B(Y) (J_nDYe-DYe,J_nDYe-DYe
)
\\
&&{} +D^2B(Y) (DYe,J_nDYe-DYe )
\\
&&{}+D^2B(Y) (J_nDYe-DYe,J_nDYe-DYe
)
\end{eqnarray*}
and these last quantities are shown to go to zero pointwise in $r$
thanks to Assumption \ref{assdphi} and to the $\alpha$-H\"
olderianity of $D^2B_n$ in the same way as for $\mC_{22}$. By
dominated convergence, $\mD_{11}$ is thus shown to converge to $0$. To
finish studying $\mD_1$ (hence $\mD$), we need to rewrite the
integrand in $\mD_{12}$ as
\begin{eqnarray*}
&& DB_n (Y_n)D^2Y_n(e,e)-DB(Y)D^2Y(e,e)
\\
&&\qquad =
DB_n(Y_n) \bigl[D^2Y_n-D^2Y
\bigr](e,e)+ \bigl[DB_n(Y_n)-DB_n(Y)
\bigr]D^2Y(e,e)
\\
&&\quad\qquad{} + \bigl[DB_n(Y)-DB(Y) \bigr]D^2Y(e,e)
\\
&&\qquad =DB_n(Y_n) \bigl[D^2Y_n-D^2Y
\bigr](e,e)+ \bigl[DB_n(Y_n)-DB_n(Y)
\bigr]D^2Y(e,e)
\\
&&\quad\qquad{} +DB(J_nY) \bigl[J_nD^2Y(e,e)-D^2Y(e,e)
\bigr]
\\
&&\quad\qquad{} + \bigl[DB(J_nY)-DB(Y) \bigr]D^2Y(e,e).
\end{eqnarray*}
The second term in this last sum is bounded in norm by
\[
\bigl\llVert D^2B_n(r,\cdot)\bigr\rrVert\llVert
Y_n-Y\rrVert\bigl\llVert D^2Y(e,e)\bigr\rrVert
\]
which goes to zero since $Y_n\to Y$ and $\llVert DB_n\rrVert $ are uniformly
bounded (as already noticed before); the norm of the third term goes to
zero because it is bounded by
\[
\bigl\llVert DB(J_nY)\bigr\rrVert\bigl\llVert J_nD^2Y(e,e)-D^2Y(e,e)
\bigr\rrVert;
\]
the norm of last term goes to zero as well by the Lipschitz property of
$DB$. Taking into account all these observations and the fact that $\mD
_{11}$ has already been shown to converge to zero, we can use
Gronwall's lemma in (\ref{eqD2}) to obtain that
\[
D^2Y_n^{s,y}(T) (e,e)-D^2Y^{s,y}(T)
(e,e)\to0.
\]
This together with the uniform boundedness of $D\Phi_n(Y_n)$ finally
yields the convergence to zero of $\mD$, hence that of the
second-order term.

At last, an application of the dominated convergence theorem with
respect to the variable $s$ in all integral terms appearing in the
Kolmogorov equation completes the proof.\quad\qed
\end{longlist}\noqed
\end{pf}
\begin{remark}
\label{remlipschitz}
Since $u$ is given as an integral of functions which are bounded in the
variable $t$, it is a Lipschitz function, hence differentiable almost
everywhere thanks to a classic result by Rademacher. Therefore, a
posteriori it satisfies the differential form of Kolmogorov equation
%
\begin{equation}
\label{eqKolmDiff} %
\qquad \cases{ \displaystyle\frac{\partial u}{\partial t}(t,y)+\bigl
\langle Du(t,y),Ay+B(t,y)\bigr\rangle+\frac{1}{2}\sum
_{j=1}^d\sigma^2_jD^2u(t,y)
(e_j,e_j)=0,
\cr
u(T,\cdot)=\Phi}
\end{equation}
for almost every $t\in[0,T]$.
\end{remark}
%
\section{Examples}
\label{secex}
We give here some examples, recalling also those mentioned at the
beginning of the paper, to which the theory exposed so far can be
applied. In particular, we first consider integral functions and
explain why they cannot be treated in the standard Hilbert space
setting for our purposes, and then show that the technical Assumption
\ref{assdphi}, which can seem very restrictive when considered in its
abstract form, is indeed satisfied by all the usual examples.

\subsection{Examples for the $\sL^p$ theory}

\begin{example}[(A negative example)]\label{ex1long}
First, we show that, as said before, even the simplest path-dependent
functions one can think of, namely integral functional, do not have
enough smoothness when considered in the standard $\sL^2$ setting.

In dimension $d=1$, consider the integral functional
\[
b_t (\gamma_t )=\int_0^tg
\bigl(\gamma(s) \bigr)\udt s,
\]
where $g:\bR\to\bR$ is a $C^3_b$ function. Its infinite-dimensional
lifting is given by
\[
B \biggl(t, \pmatrix{x
\cr
\phi} \biggr)= \pmatrix{\hat b \biggl(t,\pmatrix{x
\cr
\phi} \biggr)
\vspace*{3pt}\cr
0},
\]
where
\[
\hat b \biggl(t,\pmatrix{x
\cr
\phi} \biggr)=\int_0^tg
\bigl(\phi(s-t) \bigr)\udt s.
\]
The second G\^{a}teaux derivative of $B$ with respect to $y={x
\choose \phi}$ is simply
\[
D^2_GB (t,y ) \biggl(\pmatrix{x_1
\cr
\psi},
\pmatrix{x_2
\cr
\chi} \biggr) =\pmatrix{\displaystyle\int_0^tg^{\prime\prime}
\bigl(\phi(s-t) \bigr)\psi(s-t)\chi(s-t)\udt s
\vspace*{3pt}\cr
0}.
\]
Given ${x_1 \choose  \psi}$, ${x_2 \choose \chi}$, it is
easy to check,
by Lebesgue theorem, that this G\^{a}teaux derivative is continuous in
$y$ in the $\mathcal{L}^{2}$ topology; with some additional effort it
can be also shown that it is uniformly continuous, in $y\in\mathcal
{L}^{2}$. Presumably, thanks to this result on $B$, with due effort it
can be shown that uniform continuity of G\^{a}teaux derivatives holds
true also for the solution $Y$ of the SDE and then for $u (
t,y ) $. However, with only such knowledge about the space
regularity of $u$, we do not know how to prove that $u$ satisfies the
Kolmogorov equation (we do not know how to control the remainders in
Taylor developments). Coherently, with the present literature on the
subject, we are able to complete the proof that $u ( t,y ) $
fulfills the Kolmogorov equation only when the second-order Fr\'{e}chet
differential is uniformly continuous [not only the G\^{a}teaux
derivative for given ${x_1 \choose \psi}$, ${x_2 \choose
\chi}$]. This
is false for $B$ as above: integral functionals are not even twice
differentiable in Fr\'echet sense in general. In order for
$D^2_GB(t,y)$ to be the second-order Fr\'echet differential of $B$ we
would need that
\[
\lim_{{\llVert w\rrVert }_{\sL^2}\to0}\frac{1}{\llVert w\rrVert
}_{\sL
^2}\bigl\llVert
DB(t,y+w) z-DB(t,y)z-D_G^2B(t,y) (z,w)\bigr\rrVert
_{\sL^2}=0
\]
uniformly in $z\in\sL^2$, that is, for $y={z \choose \phi}$,
$z= {x_1 \choose \psi}$, $w={x_2 \choose \chi}$,
\begin{eqnarray*}
&&\lim_{{\llVert \chi\rrVert }_{L^2}\to0}\frac{1}{\llVert \chi\rrVert
}_{L^2}\biggl\llvert
\int_0^t \bigl[g^\prime\bigl(
\phi(s-t) +\chi(s-t) \bigr)-g^\prime\bigl(\phi(s-t) \bigr) \bigr]
\psi(s-t)\udt s
\\
&&\qquad{}-\int_0^tg^{\prime\prime} \bigl(
\phi(s-t) \bigr)\psi(s-t)\chi(s-t)\udt s\biggr\rrvert=0
\end{eqnarray*}
uniformly in $\psi\in L^2$.
Suppose that $g^{\prime\prime}$ is not constant, take as $\phi$ any
continuous function and choose $\psi(s)=s^{-\sfrac{1}{3}}$ and
$\chi_n(s)=s^{-\sfrac{1}{3}}\ind_{[-\sfrac{1}{n},0)}(s)$.
Then $\chi_n\to0$ in $L^2$ as $n\to\infty$ and
\begin{eqnarray*}
&&\lim_{n\to\infty}\frac{1}{\llVert \chi_n\rrVert }_{L^2}\biggl\llvert
\int_0^t \bigl[g^\prime\bigl(
\phi(s-t)+\chi_n(s-t) \bigr)-g^\prime\bigl(\phi(s-t) \bigr)
\bigr]\psi(s-t)\udt s
\\
&&\quad{}-\int_0^tg^{\prime\prime}
\bigl(\phi(s-t) \bigr)\psi(s-t)\chi_n(s-t)\udt s\biggr\rrvert
\\
&&\qquad =\lim_{n\to\infty}\frac{1}{\llVert \chi_n\rrVert }_{L^2}\biggl\llvert
\int_0^t \biggl[g^{\prime\prime} \bigl(
\phi(s-t) \bigr)\chi_n(s-t)\psi(s-t)
\\
&&\quad\qquad{}+\frac{1}{2}g^{\prime\prime\prime}(
\bar x)\chi_n(s-t)^2\psi(s-t) \biggr]\udt s
\\
&&\quad\qquad{}-\int
_0^tg^{\prime\prime} \bigl(\phi(s-t) \bigr)
\psi(s-t)\chi_n(s-t)\udt s\biggr\rrvert,
\end{eqnarray*}
where $\bar x$ is some point in $\bR$. Since $g^{\prime\prime\prime
}$ is bounded, we have to compute
\[
\lim_{n\to\infty}\frac{1}{\llVert \chi_n\rrVert }_{L^2}\int
_0^t\bigl\llvert\chi_n(s-t)\bigr
\rrvert^2\bigl\llvert\psi(s-t)\bigr\rrvert\udt s,
\]
but with our choice of $\chi_n$ and $\psi$ the functions $\llvert
\chi_n\rrvert ^2\llvert \psi\rrvert $ are not integrable for any~$n$. Therefore, $D^2_GB(t,y)$ cannot be the differential of second
order of $B$ in Fr\'echet sense.
\end{example}

\begin{example}
\label{ex2long}
On the other hand, the infinite-dimensional lifting of integral
functionals of the form
\[
b_{t} (\gamma_{t} ) =\int_{0}^{t}g
\bigl(\gamma(t ),\gamma(s ) \bigr)\udt s
\]
with $g$ of class $C^{2,\alpha}_b (\bR^d\times\bR^d;\bR
)$ satisfy the assumptions of Theorem \ref{thmLp} for $p=2+\alpha$;
in particular they are twice Fr\'{e}chet differentiable with
$\alpha$-H\"{o}lder continuous (hence uniformly continuous) second Fr\'{e}chet
differential in $\mathcal{L}^{p}$ for $p=2+\alpha$. Indeed, for
$y={x \choose \phi} $,
\begin{eqnarray*}
&\displaystyle B (t,y )= \pmatrix{ \displaystyle\int_{0}^{t}g \bigl(x,
\varphi(s-t ) \bigr)\udt s
\vspace*{3pt}\cr
0},&
\\
&\displaystyle D^{2}B (t,y ) \biggl(\pmatrix{x_1
\cr
\psi},
\pmatrix{x_2
\cr
\chi} \biggr)=\pmatrix{a
\cr
0},&
\end{eqnarray*}
where (denoting by $\partial_{1}$ and $\partial_{2}$ the partial derivatives
of $g$ in its two arguments)
\begin{eqnarray*}
a & =&\int_{0}^{t}\partial_{1}^{2}g
\bigl( x,\varphi(s-t ) \bigr)\udt s+\int_{0}^{t}
\partial_{2}^{2}g \bigl( \varphi(s-t ) \bigr) \psi(s-t )
\chi(s-t )\udt s
\\
&&{} +\int_{0}^{t}
\partial_{1}\partial_{2}g \bigl(x,\varphi(s-t ) \bigr)
\bigl(\psi(s-t ) +\chi(s-t ) \bigr)\udt s.
\end{eqnarray*}
For $z={x_1 \choose \phi_1}$, we have to estimate $\llVert
D^{2}B ( t,y )-D^{2}B ( t,z ) \rrVert
_{L ( \mathcal{L}^{p},\mathcal{L}^{p};\mathcal{L}^{p} )
}$ and the most difficult term is
\begin{eqnarray*}
&& \biggl\llvert\int_{0}^{t} \bigl(
\partial_{2}^{2}g \bigl( \varphi(s-t ) \bigr) -
\partial_{2}^{2}g \bigl( \varphi_{1} ( s-t )
\bigr) \bigr) \psi( s-t ) \chi( s-t )\udt s\biggr\rrvert
\\
&&\qquad \leq\bigl
\llVert
\partial_{2}^{2}g\bigr\rrVert_{\alpha}\int
_{0}^{t}\bigl\llvert\varphi( s-t ) -
\varphi_{1} ( s-t )\bigr\rrvert^{\alpha}\bigl\llvert\psi( s-t
) \bigr\rrvert\bigl\llvert\chi( s-t ) \bigr\rrvert\udt s
\\
&&\qquad \leq\bigl
\llVert
\partial_{2}^{2}g\bigr\rrVert_{\alpha}\bigl\llVert
\llvert\varphi-\varphi_{1}\rrvert^{\alpha}\bigr\rrVert
_{L^{p/\alpha}}\llVert\psi\rrVert_{L^{p}}\llVert\chi\rrVert
_{L^{p}}
\\
&&\qquad =\bigl\llVert\partial_{2}^{2}g\bigr
\rrVert_{\alpha}\llVert\varphi-\varphi_{1}\rrVert
_{L^{p}}^{\alpha}\llVert\psi\rrVert_{L^{p}}\llVert\chi
\rrVert_{L^{p}}
\end{eqnarray*}
which implies
\begin{eqnarray*}
&& \mathop{\sup_{\chi,\psi\in\mathcal{L}^{p}}}_{\llVert \chi
\rrVert _{L^{p}},\llVert \psi\rrVert _{L^{p}}\leq
1} \biggl\llvert\int
_{0}^{t} \bigl( \partial_{2}^{2}g
\bigl( \varphi( s-t ) \bigr)-\partial_{2}^{2}g \bigl(
\varphi_{1} ( s-t ) \bigr) \bigr)\psi( s-t ) \chi( s-t )\udt s\biggr
\rrvert
\\
&&\qquad  \leq\bigl\llVert\partial_{2}^{2}g\bigr\rrVert
_{C^{\alpha}}\llVert\varphi-\varphi_{1}\rrVert
_{L^{p}}^{\alpha}.
\end{eqnarray*}
Since $g$ and its derivatives are bounded, Assumption \ref{assB} is
easily seen to be satisfied.

This argument can be easily extended to include dependence on $t$ and
$s$ in $g$, as in example (i) in the \hyperref[secintro]{Introduction}.
\end{example}
%
\subsection{Examples for the theory in $\sD$}

\begin{example}
\label{ex3long}
We show now that the lifting of the function introduced in Section~\ref
{subsecmainresults}, example (ii) satisfies the assumptions of
Theorem \ref{thmmain}. For simplicity, we evaluate any c\`adl\`ag
curve $\gamma$ only in two fixed points $t_1$ and $t_2$, $0\leq
t_1\leq t_2<T$, that is, we set
\[
b_t (\gamma_t )=h_1 \bigl(
\gamma(t_1) \bigr)\ind_{[t_1,t_2)}(t)+h_2 \bigl(
\gamma(t_1),\gamma(t_2) \bigr)\ind_{[t_2,T]}(t),
\]
where $h_1:\bR^d\rightarrow\bR^d$ and $h_2:\bR^d\times\bR
^d\rightarrow\bR^d$ are in $C^{2,\alpha}_b$ on their respective
domains.

Given an element ${x \choose \phi}\in\sD$, we will write $\phi(0)$
for $x$ to avoid the burdensome notation $\phi(s)\ind
_{[-T,0)}(s)+x\ind_{\{0\}}(s)$ in the\vspace*{2pt} following computations, and we
will write $\ind_{[a,0]}$ for ${1 \choose \ind_{[a,0)}}$
accordingly.

We first check that Assumption \ref{assdphi} is satisfied. Here,
$\hat b$ is given by
\begin{eqnarray*}
\hat b_t (t,x,\phi)&=&h_1 \bigl(\phi(t_1-t
) \bigr)\ind_{[t_1,t_2)}(t)
 +h_2 \bigl(\phi(t_1-t ),\phi
(t_2-t ) \bigr)\ind_{[t_2,T]}(t).
\end{eqnarray*}
Therefore, the Fr\'echet differential of $B$ with respect to its second
argument ${x \choose \phi}$ is given by
\[
DB \biggl(t,\pmatrix{x
\cr
\phi} \biggr) \pmatrix{x_1
\cr
\psi}= \pmatrix{ D\hat b
\biggl(t,\pmatrix{x
\cr
\phi} \biggr)\pmatrix{x_1
\cr
\psi}
\vspace*{3pt}\cr
0},
\]
where
\begin{eqnarray*}
&& D\hat b \biggl(t,\pmatrix{x
\cr
\phi} \biggr)\pmatrix{x_1
\cr
\psi}
\\
&&\qquad = D
h_1 \bigl(\phi(t_1-t ) \bigr)\psi(t_1-t )
\ind_{[t_1,t_2)}(t)
\\
&&\quad\qquad{}+D h_2 \bigl(\phi(t_1-t ),\phi
(t_2-t ) \bigr) \bigl(\psi(t_1-t ),\psi
(t_2-t ) \bigr)\ind_{[t_2,T]}(t)
\end{eqnarray*}
and $Dh_j$ denotes the Jacobian matrix of $h_j$.

For any fixed $a\in[-T,0]$ (recall the convention we adopted in
Definition \ref{defphi}), the first component of $DB (t,
{x \choose \phi} )J_n {1 \choose \ind_{[a,0)}}$ is given by
\begin{eqnarray*}
&& \bigl[D h_1 \bigl(\phi(t_1-t ) \bigr)\cdot
J_n\ind_{[a,0]} (t_1-t ) \bigr]
\ind_{[t_1,t_2)}(t)
\\
&&\qquad{} + \bigl[D h_2 \bigl(\phi(t_1-t ),\phi
(t_2-t ) \bigr)\cdot\bigl(J_n\ind_{[a,0]}
(t_1-t ),J_n\ind_{[a,0]} (t_2-t )
\bigr) \bigr]\ind_{[t_2,T]}(t)
\end{eqnarray*}
while the second is $0$. Therefore,
\begin{eqnarray*}
&&DB \biggl(t,\pmatrix{x
\cr
\phi} \biggr) \biggl(J_n\pmatrix{1
\cr
\ind
_{[a,0)}}-\pmatrix{1
\cr
\ind_{[a,0)}} \biggr)\longrightarrow0
\end{eqnarray*}
if and only if
\begin{eqnarray*}
&&J_n\ind_{[a,0]} (t_j-t )\rightarrow
\ind_{[a,0]} (t_j-t ),
\end{eqnarray*}
$j=1,2$. Fix $j=1$ (the situation being analogous with $j=2$); if
$t=t_1$, it is straightforward to verify the assumption, therefore,
suppose $t\neq t_1$. Then, using the sequence $J_n$ given by (\ref
{eqJn}), if $t_1>0$ we have
%
\begin{eqnarray}
\label{eqJnex2} J_n\ind_{[a,0]} (t_1-t )&=&\int
_{-T}^0\rho_n \bigl(\tau
_{\sfrac{1}{n}} (t_1-t )-y \bigr)\ind_{[a,0]}(y)\udt y
\nonumber\\[-8pt]\\[-8pt]
\nonumber
&=&\int_a^0\rho_n
(t_1-t-y )\udt y
\end{eqnarray}
for $n$ big enough. Now if $t_1-t<a$ then choosing $n$ large enough we
have that $ (t_1-t )+\supp(\rho_n )\cap
[a,0]=\varnothing$, hence the function in (\ref{eqJnex2}) equals to $0$
definitively as $n$ tends to infinity. Conversely, if $t_1-t>a$ for $n$
large enough we have that $ (t_1-t )+\supp(\rho
_n )\cap[a,0]= (t_1-t )+\supp(\rho_n )$
and the function in (\ref{eqJnex2}) equals $1$ definitively. If
$t_1=0$, the same procedure applies when $t\neq T$ or $a>-T$, while
when $t=T$ and $a=-T$ by the definition of $\tau_{\sfrac{1}{n}}$ it
follows that
\begin{eqnarray*}
J_n\ind_{[a,0]} (-T )&=&\int_{-T}^0
\rho_n \bigl(\tau_{\sfrac{1}{n}} (-T )-y \bigr)\ind_{[-T,0)}(y)\udt y
\\
&=&\int_{-T}^0\rho_n \biggl(-T+
\frac{1}{n}-y \biggr)\udt y=1.
\end{eqnarray*}
Therefore, for any $t\in[0,T]$, for any $a\neq t_1-t$ we have that
$J_n\ind_{[a,0]} (t_1-t )=\ind_{[a,0]} (t_1-t
)$ definitively as $n$ tends to $\infty$, as required. It is easy to
see that if $a=t_1-t$ then $J_n\ind_{[a,0]} (t_1-t )\to
\frac{1}{2}$.

The second Fr\'echet differential is given by
\[
D^2B \biggl(t,\pmatrix{x
\cr
\phi} \biggr) \biggl(\pmatrix{x_1
\cr
\psi} ,
\pmatrix{x_2
\cr
\chi} \biggr)= \pmatrix{D^2\hat b \biggl(t,\pmatrix{x
\cr
\phi} \biggr) \biggl(\pmatrix{x_1
\cr
\psi} , \pmatrix{x_2
\cr
\chi}\biggr)
\vspace*{3pt}\cr
0},
\]
where
\begin{eqnarray*}
&&D^2\hat b \biggl(t,\pmatrix{x
\cr
\phi} \biggr) \biggl(\pmatrix{x_1
\cr
\psi},\pmatrix{x_2
\cr
\chi} \biggr)
\\
&&\qquad =D^2h_1 \bigl(
\phi(t_1-t) \bigr) \bigl(\psi(t_1-t),
\chi(t_1-t) \bigr)\ind_{[t_1,t_2)}(t)
\\
&&\quad\qquad{}+D^2h_2 \bigl(\phi(t_1-t),
\phi(t_2-t) \bigr)
\\
&&\quad\qquad{}\times \bigl( \bigl(\psi(t_1-t),
\psi(t_2-t) \bigr), \bigl(\chi(t_1-t),
\chi(t_2-t) \bigr) \bigr)
\ind_{[t_2,T]}(t)
\end{eqnarray*}
and $D^2h_j$ denotes the Hessian tensor of $h_j$; it can be easily seen
that this differential satisfies the requirements of Assumption \ref
{assdphi} reasoning as above.

It is also immediate to check that since $h_1$ and $h_2$ are in
$C^{2,\alpha}_b$ Assumption \ref{assB} is satisfied by this example.
\end{example}

\begin{example}
\label{ex4long}
We can use evaluation at fixed times also in the terminal condition for
the path-dependent Kolmogorov equation (\ref{eqPKolmogorovintroCF})
(see also Section~\ref{seccomparison}): given a smooth function
$q:\bR^{(n+1)d}\to\bR$, bounded with bounded derivatives, consider
\[
f (\gamma_T )=q \bigl(\gamma(t_0),\gamma(t_1),
\dots,\gamma(t_n),\gamma(T) \bigr).
\]
Its infinite-dimensional lifting is then given by
\[
\Phi\pmatrix{x
\cr
\phi}= \pmatrix{\hat f\pmatrix{x
\cr
\phi}
\vspace*{3pt}\cr
0},
\]
where
\[
\hat f \biggl(\pmatrix{x
\cr
\phi} \biggr)=q \bigl(\phi(t_0-T ),\phi
(t_1-T ),\dots,\phi(t_n-T ),x \bigr).
\]
From Example \ref{ex3long}, it is immediate to see that such a $\Phi
$ satisfies Assumption \ref{assdphi} and, therefore, it can be chosen
as terminal condition in Theorem \ref{thmmain}.
\end{example}

\begin{example}
From Examples \ref{ex3long} and \ref{ex4long}, it follows also that
Theorem \ref{thmmain} can be applied when the drift or the terminal
condition in the Kolmogorov equation (or both) are delayed functions of
the form
\[
b_t (\gamma_t )=g \bigl(\gamma(t),\gamma(t-\delta)
\bigr)\ind_{[\delta,T]}(t), \qquad f (\gamma_T )=q \bigl(\gamma(T),
\gamma(T-\delta) \bigr)
\]
for $g$ and $q$ sufficiently regular and with values in $\bR^d$ and
$\bR$, respectively, and $0<\delta<T$, since in this case we have that
\[
B \biggl(t,\pmatrix{x
\cr
\phi} \biggr)= \pmatrix{g \bigl(x,\phi(-\delta) \bigr)
\cr
0}
\ind_{[\delta,T]}(t)\qquad \forall t\in[0,T]
\]
and
\[
\Phi\pmatrix{x
\cr
\phi}=\pmatrix{q \bigl(x,\phi(-\delta) \bigr)
\cr
0}.
\]
\end{example}

\begin{remark}
The theory exposed here cannot be applied to example (iv) in
Section~\ref{subsecmainresults}, that is the functional
\[
b_t (\gamma_t )=\sup_{s\in[0,t]}\gamma(s)
\]
since the supremum is not Fr\'echet differentiable as a function of the path.
\end{remark}
%
\section{Comparison with path-dependent calculus}\label{seccomparison}
We conclude this work establishing some connections between our results
and objects and those defined by Dupire and successively developed by
Cont and Fourni\'e.
We recall here the definitions of the pathwise derivatives given in
\citet{CF1}. For a function $\nu= \{\nu_t \}_t$,
$\nu_t\colon D ([0,T];\bR^d )\to\bR^d$, the $i$th \emph
{vertical derivative} at $\gamma_t$ ($i=1,\dots,d$) is defined as
%
\begin{equation}
\label{eqdefvertD} \rD_i\nu_t(\gamma_t)=\lim
_{h\to0}\frac{\nu_t (\gamma
_t^{he_i} )-\nu_t(\gamma_t)}{h},
\end{equation}
where $\gamma_t^{he_i}(s)=\gamma_t(s)+he_i\ind_{\{t\}}(s)$; we
denote the \emph{vertical gradient} at $\gamma_t$ by
\[
\rD\nu_t(\gamma_t)= \bigl(\rD_1
\nu_t(\gamma_t),\dots,\rD_d\nu_t(
\gamma_t) \bigr);
\]
higher order vertical derivatives are defined in a straightforward way.
The \emph{horizontal derivative} at $\gamma_t$ is defined as
%
\begin{equation}
\label{eqdefhorizD} \rD_t\nu(\gamma_t )=\lim
_{h\to0^+}\frac{\nu
_{t+h} (\gamma_{t,h} )-\nu_t (\gamma_t )}{h},
\end{equation}
where $\gamma_{t,h}(s)=\gamma_t(s)\ind_{[0,t]}(s)+\gamma_t(t)\ind
_{(t,t+h]}(s)\in D ([0,t+h];\bR^d )$.
The connection between a functional $b$ of paths and the operator $B$
was essentially a matter of definition, as carried out in (\ref
{eqMt})--(\ref{eqb}). To establish some relations between Fr\'echet
differentials of $B$ and horizontal and vertical derivatives of $b$ is
much less obvious; some results are given by the following theorem.
\begin{theorem}
\label{thmkolm}
Suppose $u\colon[0,T]\times\sD\to\bR$ is given and define, for
each $t\in[0,T]$, $\nu_t\colon D([0,t];\bR^d)\to\bR$ as $\nu
_t(\gamma):=u(t,\gamma(t),L_t\gamma)$, in the same way as in~(\ref
{eqb}). Then the vertical derivatives of $\nu_t$ coincide with the
partial derivatives of $u$ with respect to the second variable (i.e.,
the present state), that is,
%
\begin{equation}
\label{eqequivhoriz} \rD_i\nu_t(\gamma)=\frac{\partial}{\partial x}u(t,x,L_t
\gamma),\qquad i=1,\dots, d.
\end{equation}
The same result holds true also if $u$ is given from $\nu$ as in (\ref{eqB}).
Furthermore let $\gamma_t\in C^1_b([0,t];\bR^d)$ and let again $u$ be
given and define $\nu$ as above. Then
\[
\rD_t\nu(\gamma_t)=\frac{\partial u}{\partial t} \bigl(t,\gamma
(t),L_t\gamma_t \bigr)+{\bigl\langle Du\bigl(t,
\gamma(t),L_t\gamma_t\bigr),(L_t
\gamma_t)^\prime_+\bigr\rangle},
\]
where $\langle\cdot,\cdot\rangle$ is the duality between $D$ and
$D^{\prime}$, $Du$ is the Fr\'echet derivative of $u$ with respect to
$\phi$ and the lower script $+$ denotes right derivative.
\end{theorem}
\begin{pf}
Both claims in the theorem are proved through explicit calculations
starting from the definition of derivatives.
From the definition of vertical derivative, one gets
\begin{eqnarray*}
\rD_i\nu_t(\gamma)&=&\lim_{h\to0}
\frac{1}{h} \bigl[\nu_t \bigl(\gamma^{h} \bigr)-
\nu_t(\gamma) \bigr]
\\
&=&\lim_{h\to0}\frac{1}{h} \bigl[u \bigl(t,
\gamma^{h}(t),L_t\gamma^{h} \bigr)-u\bigl(t,
\gamma(t),L_t\gamma\bigr) \bigr]
\\
&=&\lim_{h\to0}\frac{1}{h} \bigl[u \bigl(t,
\gamma(t)+h,L_t\gamma^{h} \bigr)-u\bigl(t,
\gamma(t),L_t\gamma\bigr) \bigr]
\\
&=&\frac{\partial}{\partial x_i}u(t,x,L_t\gamma).
\end{eqnarray*}
This proves the first part of the theorem.

For the second part, suppose first that there is no explicit dependence
on $t$ in $u$. Then
\begin{eqnarray*}
\rD_t\nu(\gamma_t)&=&\lim_{h\to0}
\frac{1}{h} \bigl[u \bigl(\gamma_{t,h}(t),L_{t+h}
\gamma_{t,h} \bigr)-u \bigl(t,\gamma_t(t),L_t
\gamma_t \bigr) \bigr]
\\
&=&\lim_{h\to0}\frac{1}{h} \bigl[u \bigl(
\gamma_t(t),L_{t+h}\gamma_{t,h} \bigr)-u \bigl(t,
\gamma_t(t),L_t\gamma_t \bigr) \bigr]
\\
&=&\lim_{h\to0}\frac{1}{h} \biggl[u
\biggl(\gamma_t(t),
\cases{
\gamma_{t,h}(t+s), &\quad $[-t-h,0)$,
\cr
\gamma_{t,h}(0),&\quad $[-T,-t-h)$}
\biggr)
\\
&&{}
-u \biggl(\gamma_t(t),
\cases{
\gamma_t(t+s), &\quad $[-t,0)$,
\cr
\gamma_t(0), &\quad $[-T,-t)$}\biggr) \biggr]
\\
&=&\lim_{h\to0}\frac{1}{h}
\left[u
\left(\gamma_t(t),
\cases{
\gamma_t(t),& \quad $[-h,0)$,
\cr
\gamma_t(t+s+h),&\quad $[-t,-h)$,
\cr
\gamma_t(t+s+h),&\quad $[-t-h,-t)$,
\cr
\gamma_t(0),&\quad $[-T,-t-h)$}
\right)\right.
\\
&&{}
-u
\left.\left(
\gamma_t(t),
\cases{
\gamma_t(t+s),&\quad $[-h,0)$,
\cr
\gamma_t(t+s),&\quad $[-t,-h)$,
\cr
\gamma_t(0),&\quad $[-t-h,-t)$,
\cr
\gamma_t(0),&\quad $[-T,-t-h)$}\right)\right].
\end{eqnarray*}
Last line can be written as
%
\begin{equation}
\label{eqsemplice} \lim_{h\to0}\frac{1}{h} \bigl[u \bigl(
\gamma_t(t),L_t\gamma_t+N_{t,h}
\gamma_t \bigr)-u \bigl(\gamma_t(t),L_t
\gamma_t \bigr) \bigr],
\end{equation}
where
%
\begin{equation}
\label{eqdefN} N_{t,h}\gamma_t(s)=
\cases{0, &\quad$[-T,-t-h)$,
\cr
\gamma_t(t+h+s)-\gamma_t(0), &\quad$[-t-h,-t)$,
\cr
\gamma_t(t+h+s)-\gamma_t(t+s), &\quad$[-t,-h)$,
\cr
\gamma_t(t)-\gamma(t+s), &\quad$[-h,0)$.}
\end{equation}
$N_{t,h}\gamma_t$ is a continuous function that goes to $0$ as $h\to
0$; moreover, recalling that in the definition of horizontal derivative
$h$ is greater than zero, we see that:
\begin{longlist}[(iii)]
\item[(i)]  for $s\in[-T,-t)$ $\exists\bar h$ s.t.
$s<-t-\bar h$, hence $N_{t,h}\gamma(s)=0$ $\forall h<\bar h$ and
\[
\lim
_{h\to0^+}\frac{1}{h}N_{t,h}\gamma(s)=0=(L_t\gamma)^\prime(s);
\]

\item[(ii)]  for $s=-t$, since $N_{t,h}\gamma(-t)=\gamma
(h)-\gamma(0)$ we have
\[
\lim
_{h\to0^+}\frac{1}{h}N_{t,h}\gamma
_t(-t)= \biggl(\frac{\ud^+}{\ud s}L_t\gamma_t
\biggr)(-t)=(L_t\gamma_t)^\prime_+(-t)=\gamma^\prime_+(0);
\]

\item[(iii)] for $s\in(-t,0)$ $\exists\bar h$ s.t. $s<-\bar
h<0$, hence
\begin{eqnarray*}
\lim_{h\to0^+}\frac{1}{h}N_{t,h}\gamma
_t(s)&=&\lim_{h\to0^+}\frac{1}{h} \bigl[\gamma_t(t+s+h)-\gamma_t(t+s)
\bigr]
\\
&=& \gamma^\prime_+(t+s)=\gamma^\prime(t+s)=(L_t\gamma
_t)^\prime(s).
\end{eqnarray*}
\end{longlist}
Therefore,
%
\begin{equation}
\label{eqconvderiv} \frac{1}{h}N_{t,h}\gamma_t(s)
\stackrel{h\to0^+} {\longrightarrow}(L_t\gamma_t)^\prime_+(s)
\end{equation}
and, since $\gamma\in C^1_b$,
\[
(L_t\gamma_t)^\prime_+(s)=(L_t
\gamma_t)^\prime(s)\qquad\forall s\neq-t.
\]
Again since $\gamma_t\in C^1$ with bounded derivative, $\frac
{1}{h}N_{t,h}\gamma_t$ converges to $(L_t\gamma_t)^\prime_+$ also uniformly.
Keeping into account (\ref{eqsemplice}) and the definition of Fr\'
echet derivative, one gets
\begin{eqnarray*}
\rD_t\nu(\gamma_t)&=&\lim_{h\to0}
\frac{1}{h} \bigl[u\bigl(\gamma_t(t),L_t
\gamma_t+N_{t,h}\gamma_t\bigr)-u\bigl(
\gamma_t(t),L_t\gamma_t\bigr) \bigr]
\\
&=&\lim_{h\to0}\frac{1}{h} \bigl[ \bigl\langle D u\bigl(
\gamma_t(t),L_t\gamma_t\bigr),
N_{t,h}\gamma_t \bigr\rangle+\xi(h) \bigr],
\end{eqnarray*}
where $\xi$ is infinitesimal with respect to ${\llVert
N_{t,h}\gamma_t\rrVert }$ as $h\to0$,
\begin{eqnarray*}
&=&\lim_{h\to0}\frac{1}{h} \bigl\langle D
u\bigl(\gamma_t(t),L_t\gamma_t
\bigr),N_{t,h}\gamma_t \bigr\rangle+\lim
_{h\to0}\frac{{\llVert N_{t,h}\gamma_t\rrVert }}{h}\frac{\xi
(h)}{{\llVert
N_{t,h}\gamma_t\rrVert }}
\\
&=&\bigl\langle D u\bigl(\gamma_t(t),L_t
\gamma_t\bigr), (L_t\gamma_t)^\prime
_+\bigr\rangle
\end{eqnarray*}
by the dominated convergence theorem.

If now $u$ depends explicitly on $t$ just write
\begin{eqnarray*}
\frac{1}{h} \bigl[\nu_{t+h}(\gamma_{t,h})-
\nu_t(\gamma) \bigr]&=&\frac{1}{h} \bigl[u \bigl(t+h,
\gamma(t),L_{t+h}\gamma_{t,h} \bigr)-u \bigl(t,
\gamma(t),L_t\gamma\bigr) \bigr]
\\
&=&\frac{1}{h} \bigl[u \bigl(t+h,\gamma(t),L_{t+h}
\gamma_{t,h} \bigr)-u \bigl(t,\gamma(t),L_{t+h}
\gamma_{t,h} \bigr) \bigr]
\\
&&{}+\frac{1}{h} \bigl[u \bigl(t,\gamma(t),L_{t+h}
\gamma_{t,h} \bigr)-u \bigl(t,\gamma(t),L_t\gamma\bigr)
\bigr];
\end{eqnarray*}
the first term in the last line converges to the time derivative of $u$
while the second can be treated exactly as above.
\end{pf}

Thanks to this result we can reinterpret equation (\ref{eqKolmDiff}),
which is the differential form of the infinite-dimensional Kolmogorov
equation (\ref{eqPKolmogorov}), in terms of the horizontal and
vertical derivatives introduced in the previous section.

Consider the Kolmogorov equation with horizontal and vertical
derivatives, namely
%
\begin{equation}
\label{eqPHVKolmogorov} %
\cases{ \displaystyle\rD_t\nu(\gamma_t)+b_t(
\gamma_t)\cdot\rD\nu_t(\gamma_t)+
\frac{1}{2}\sum_{j=1}^d
\sigma_j^2\rD_j^2
\nu_t(\gamma_t)=0,
\cr
\nu_T(
\gamma_T)=f(\gamma_T).}
\end{equation}
\begin{theorem}
\label{thmpdkolm}
Let $X^{\gamma_t}$ be the solution to equation
%
{\renewcommand{\theequation}{5}
\begin{equation}
\cases{ \ud X(t)=b_t(X_t)\udt t+\sigma\udt
W(t), &\quad for $t\in[t_0,T]$,
\cr
X_{t_0}=
\gamma_{t_0}.}
\end{equation}\setcounter{equation}{54}}%
Associate to $b_t$ and $f$ the operators $B$ and $\Phi$ as in (\ref
{eqdefB}); if such $B$ and $\Phi$ satisfy the assumptions of Theorem
\ref{thmmain} then, for almost every $t$, the function
%
\begin{equation}
\label{eqsolnu} \nu_t(\gamma_t)=\bE\bigl[f
\bigl(X^{\gamma_t}(T) \bigr) \bigr]
\end{equation}
is a solution of the path-dependent Kolmogorov equation (\ref
{eqPHVKolmogorov}) for all $\gamma\in C^1_b ([0,T];
\bR^d
)$ such that $\gamma^\prime(0)=0$.
\end{theorem}
\begin{pf}
Lift equation (\ref{eqPSDE}) to the infinite-dimensional SDE (\ref
{eqDYshort}) defining the operators $A$, $B$ and $\Sigma$ as in the
previous sections; associate then to this last equation the PDE (\ref
{eqPKolmogorov}) with terminal condition given by
\[
\Phi\biggl(\pmatrix{x
\cr
\phi} \biggr)=f \biggl(\widetilde M\pmatrix{x
\cr
\phi} \biggr).
\]
Fix $t$: with\vspace*{1pt} our choice of $\gamma$ the element $y=(\gamma
(t),L_t\gamma_t)$ is in $\Dom(A_{\aC} )$, therefore, if
$B$ and $\Phi$ satisfy Assumptions \ref{assB} and \ref{assdphi},
Theorem \ref{thmmain} guarantees that $u(s,y)=\bE[\Phi
(Y^{s,y}(T) ) ]$ is a solution to the Kolmogorov equation.
Notice that solving this equation for $s\geq t$ involves only a piece
(possibly all) of the path $\gamma_t$, so that our ``artificial''
lengthening by means of $L_t$ is used only for defining all objects in
the right way but does not come into the solution of the equation. Of
course, in principle one can solve the infinite-dimensional PDE for any
$s\in[0,T]$, anyway we are interested in solving it at time $t$:
indeed if we now define $\nu$ through $u$ by means of (\ref{eqb}) we
have that
\begin{eqnarray*}
\nu_t(\gamma_t)&=&u\bigl(t,\gamma(t),L_t
\gamma_t\bigr)
\\
&=&\bE\bigl[f \bigl(\widetilde M \bigl(Y^{t,y}(T) \bigr) \bigr) \bigr]
\\
&=&\bE\bigl[f \bigl(X^{\gamma_t}(T) \bigr) \bigr].
\end{eqnarray*}
Recalling Remark \ref{remlipschitz} and noticing that $
(L_t\gamma_t )^\prime_+=A (L_t\gamma_t )$ thanks to
the assumption that $\gamma^\prime(0)=0$, we can apply for almost
every $t$ Theorem \ref{thmkolm} obtaining that equations (\ref
{eqKolmDiff}) and (\ref{eqPHVKolmogorov}) coincide.
\end{pf}
\begin{remark}
If in the above proof one can show that the function $u$ which solves
(\ref{eqPKolmogorov}) is in fact differentiable with respect to $t$
for every $t\in[0,T]$, then Theorem \ref{thmpdkolm} holds
everywhere, that is, the function $\nu$ defined by (\ref{eqsolnu})
solves equation~(\ref{eqPHVKolmogorov}) for every $t\in[0,T]$.
\end{remark}
\begin{remark}
The restriction $\gamma^\prime(0)=0$ is only technical and is likely
avoidable with some effort. We intend to address this matter in the
future to obtain full generality in our result.
\end{remark}

\begin{appendix}\label{append}
\section*{Appendix: Proof of Theorem \texorpdfstring{\protect\ref{4535454754545}}{2.3}}
Thanks to Theorem \ref{thmSDE} we can work path by path. Therefore,
we consider $\omega$ fixed throughout the proof.

\renewcommand{\theequation}{A\arabic{equation}}
\setcounter{equation}{0}
We start from a simple estimate; for $y$, $k\in E$ we have
\begin{eqnarray*}
&& \bigl\llVert Y^{t_0,y+k}(t)-Y^{t_0,y}(t)\bigr\rrVert
{_E}
\\
&&\qquad =\biggl\llVert e^{(t-t_0)A}k+\int_{t_0}^te^{(t-s)A}
\bigl[B \bigl(s,Y^{t_0,y+k}(s) \bigr)-B \bigl(s,Y^{t_0,y}(s) \bigr)
\bigr]\udt s\biggr\rrVert_E
\\
&&\qquad \leq C{\llVert k\rrVert}_E+C{\llVert DB\rrVert}_\infty
\int_{t_0}^t{\bigl\llVert Y^{t_0,y+k}(s)-Y^{t_0,y}(s)
\bigr\rrVert}_E\udt s
\end{eqnarray*}
hence, by Gronwall's lemma,
%
\begin{equation}
\label{eqstimaY} \sup_t{\bigl\llVert
Y^{t_0,y+k}(t)-Y^{t_0,y}(t)\bigr\rrVert}_E\leq\widetilde
C_Y{\llVert k\rrVert}_E.
\end{equation}

%
\textit{First derivative}.
We introduce the following equation for the unknown $\xi^{t_0,y}(t)$
taking values in the space of linear bounded operators $L(E,E)$
\[
\xi^{t_0,y}(t)=e^{(t-t_0)A}+\int_{t_0}^te^{(t-s)A}DB
\bigl(s,Y^{t_0,y}(s) \bigr)\xi^{t_0,y}(s)\udt s.
\]
Existence and uniqueness of a solution in $L^\infty
(0,T;L(E,E) )$ follow again easily from the contraction mapping
principle, since
\begin{eqnarray*}
&& \biggl\llVert\int_{t_0}^te^{(t-s)A}DB
\bigl(s,Y^{t_0,y}(s) \bigr) \bigl[\xi_1(s)-\xi
_2(s) \bigr]\udt s\biggr\rrVert_{L(E,E)}
\\
&&\qquad \leq C{\llVert DB\rrVert}_\infty\int_{t_0}^t{
\bigl\llVert\xi_1(s)-\xi_2(s)\bigr\rrVert
}_{L(E,E)}\udt s.
\end{eqnarray*}
Moreover, by Gronwall's lemma, ${\llVert \xi^{t_0,y}(t)\rrVert
}_{L(E,E)}\leq
C_\xi$ uniformly in $t$.
Now for $k\in E$ we compute
\begin{eqnarray*}
r^{t_0,y,k}(t)&:=&Y^{t_0,y+k}(t)-Y^{t_0,y}(t)-
\xi^{t_0,y}(t)k
\\[-2pt]
&=&\int_{t_0}^te^{(t-s)A}
\bigl[B \bigl(s,Y^{t_0,y+k}(s) \bigr)-B \bigl(s,Y^{t_0,y}(s) \bigr)
\bigr]\udt s
\\[-2pt]
&&{}-\int_{t_0}^te^{(t-s)A}DB
\bigl(s,Y^{t_0,y}(s) \bigr)\xi^{t_0,y}(s)k\udt s
\\[-2pt]
&=&\int_{t_0}^te^{(t-s)A}
\biggl[\int_0^1 DB \bigl(s,\alpha
Y^{t_0,y+k}(s)
\\[-2pt]
&&{}+(1-\alpha)Y^{t_0,y}(s) \bigr) \bigl(Y^{t_0,y+k}(s)-Y^{t_0,y}(s)
\bigr)\udt\alpha
\\[-2pt]
&&{}-DB \bigl(s,Y^{t_0,y}(s) \bigr)\xi^{t_0,y}(s)k
\biggr]\udt s
\\[-2pt]
&=&\int_{t_0}^te^{(t-s)A}DB
\bigl(s,Y^{t_0,y}(s) \bigr)r^{t_0,y,k}(s)\udt s
\\[-2pt]
&&{}+\int_{t_0}^te^{(t-s)A}
\biggl[\int_0^1 DB \bigl(s,\alpha
Y^{t_0,y+k}(s)+(1-\alpha)Y^{t_0,y}(s) \bigr)\udt\alpha
\\[-2pt]
&&{}-DB \bigl(s,Y^{t_0,y}(s) \bigr) \biggr]
\bigl(Y^{t_0,y+k}(s)-Y^{t_0,y}(s) \bigr)\udt s.
\end{eqnarray*}
Recalling (\ref{eqstimaY}), we get
\begin{eqnarray*}
\bigl\llVert r^{t_0,y,k} (t)\bigr\rrVert_E
&\leq& C{\llVert DB\rrVert}_\infty\int_{t_0}^t{
\bigl\llVert r^{t_0,y,k}(s)\bigr\rrVert}_E\udt s
\\[-2pt]
&&{}+C\cdot\widetilde C_Y{\llVert k\rrVert
}_E\int_{t_0}^t\biggl\llVert\int
_0^1 DB \bigl(s,\alpha Y^{t_0,y+k}(s)+(1-
\alpha)Y^{t_0,y}(s) \bigr)\udt\alpha
\\[-2pt]
&&{}-DB \bigl(s,Y^{t_0,y}(s) \bigr)
\biggr\rrVert_{L(E,E)}\udt s
\\[-2pt]
&\leq& C{\llVert DB\rrVert}_\infty\int_{t_0}^t{
\bigl\llVert r^{t_0,y,k}(s)\bigr\rrVert}_E\udt s
\\[-2pt]
&&{}+C\cdot\widetilde C_Y{\llVert k\rrVert
}_E{\bigl\llVert D^2B\bigr\rrVert}_\infty\int
_{t_0}^t\int_0^1
\alpha{\bigl\llVert Y^{t_0,y+k}(s)+Y^{t_0,y}(s)\bigr\rrVert
}_E\udt\alpha\udt s
\\[-2pt]
&\leq& C{\llVert DB\rrVert}_\infty\int_{t_0}^t{
\bigl\llVert r^{t_0,y,k}(s)\bigr\rrVert}_E\udt s+C\cdot\widetilde
C_Y(T-t_0){\bigl\llVert D^2B\bigr\rrVert
}_\infty{\llVert k\rrVert}_E^2
\end{eqnarray*}
which yields, by Gronwall's lemma,
\[
{\bigl\llVert r^{t_0,y,k}(t)\bigr\rrVert}_E\leq\widetilde C{
\llVert k\rrVert^2}_E.
\]
Therefore,
\[
\xi^{t_0,y}(t)k=DY^{t_0,y}(t)k \qquad\forall k\in E.
\]
We proceed with an estimate about the continuity of $\xi^{t_0,y}(t)$
with respect to the initial condition $y$. For $h$, $k\in E$
\begin{eqnarray*}
&& \bigl\llVert\xi^{t_0,y+k}(t)h -\xi^{t_0,y}(t)h\bigr\rrVert
_E
\\[-2pt]
&&\qquad =\biggl\llVert\int_{t_0}^te^{(t-s)A}
\bigl[DB \bigl(s,Y^{t_0,y+k}(s) \bigr)\xi^{t_0,y+k}(s)h-DB
\bigl(s,Y^{t_0,y}(s) \bigr)\xi^{t_0,y}h \bigr]\udt s\biggr\rrVert
_E
\\[-2pt]
&&\qquad \leq\biggl\llVert\int_{t_0}^te^{(t-s)A}
\bigl[DB \bigl(s,Y^{t_0,y+k}(s) \bigr)\xi^{t_0,y+k}(s)h
\\[-2pt]
&&\quad\qquad{}-DB
\bigl(s,Y^{t_0,y+k}(s) \bigr)\xi^{t_0,y}(s)h \bigr]\udt s\biggr\rrVert
_E
\\[-2pt]
&&\quad\qquad{}+\biggl\llVert\int_{t_0}^te^{(t-s)A}
\bigl[DB \bigl(s,Y^{t_0,y+k}(s) \bigr)\xi^{t_0,y}(s)h
\\[-2pt]
&&\quad\qquad{}-DB
\bigl(s,Y^{t_0,y}(s) \bigr)\xi^{t_0,y}(s)h \bigr]\udt s\biggr\rrVert
_E
\\[-2pt]
&&\qquad \leq C{\llVert DB\rrVert}_\infty\int
_{t_0}^t\bigl\llVert\xi^{t_0,y+k}(s)h-\xi
^{t_0,y}(s)h\bigr\rrVert_E\udt s
\\[-2pt]
&&\quad\qquad{}+C\int_{t_0}^t
\bigl\llVert DB \bigl(s,Y^{t_0,y+k}(s) \bigr)-DB \bigl(s,Y^{t_0,y}(s)
\bigr)\bigr\rrVert_{L(E,E)}
\bigl
\llVert\xi^{t_0,y}(s)h\bigr\rrVert_E\udt s
\\[-2pt]
&&\qquad \leq C{\llVert DB
\rrVert}_\infty\int_{t_0}^t\bigl\llVert
\xi^{t_0,y+k}(s)h-\xi^{t_0,y}(s)h\bigr\rrVert_E\udt s
\\[-2pt]
&&\quad\qquad{}
+C\cdot C_\xi\llVert h\rrVert_E\bigl\llVert
D^2B\bigr\rrVert_\infty\int_{t_0}^t
\bigl\llVert Y^{t_0,y+k}(s)-Y^{t_0,y}(s)\bigr\rrVert
_E\udt s
\\[-2pt]
&&\qquad \leq C{\llVert DB\rrVert}_\infty\int
_{t_0}^t\bigl\llVert\xi^{t_0,y+k}(s)h-\xi
^{t_0,y}(s)h\bigr\rrVert_E\udt s
\\[-2pt]
&&\quad\qquad{}+C\cdot C_\xi{
\bigl\llVert D^2B\bigr\rrVert}_\infty\widetilde
C_Y(t-t_0){\llVert h\rrVert}_E{\llVert
k\rrVert}_E.
\end{eqnarray*}
Again by Gronwall's lemma, we get
%
\begin{equation}
\stepcounter{equation}
\label{eqstimaDxi} {\bigl\llVert
\xi^{t_0,y+k}(t)h-\xi^{t_0,y}(t)h\bigr\rrVert}_E\leq
\widetilde C_\xi{\llVert h\rrVert}_E{\llVert k\rrVert
}_E.
\end{equation}
Therefore, $\xi^{t_0,y}(t)$ is uniformly continuous in $y$ uniformly
in $t$.\vspace*{6pt}

\textit{Second derivative}.
Let us consider the operator $\sU$ defined on the space $C
([t_0,T]$; $L(E,E;E) )$ through the equation
\begin{eqnarray*}
\sU(Y) (t) (h,k)&=&\int_{t_0}^te^{(t-s)A}D^2B
\bigl(s,Y^{t_0,y}(s)\bigr) \bigl(\xi^{t_0,y}(s)h,
\xi^{t_0,y}(s)k\bigr)\udt s
\\
&&{}+\int_{t_0}^te^{(t-s)A}DB
\bigl(s,Y^{t_0,y}(s)\bigr)Y(s) (h,k)\udt s
\end{eqnarray*}
for $h$, $k\in E$; it is immediate to check that $\sU(Y)$ belongs to
$C ([t_0,T];L(E,E;E) )$.

Since
\begin{eqnarray*}
&& \sup_{t,h,k}\bigl\llVert\sU(Y_1) (t) (h,k)-
\sU(Y_2) (t) (h,k)\bigr\rrVert_E
\\
&&\qquad \leq C{\llVert DB\rrVert}_\infty T\sup_{t,h,k}\bigl
\llVert Y_1(t) (h,k)-Y_2(t) (h,k)\bigr\rrVert
_E
\end{eqnarray*}
there exists a unique fixed point for $\sU$, which will be denoted by
$\eta^{t_0,y}(t)(h,k)$; furthermore simple calculations yield that
$\llVert \eta^{t_0,y}(t)\rrVert _{L(E,E;E)}\leq C_\eta$
uniformly in $t$. We now compute:
\begin{eqnarray*}
\tilde r^{t_0,y,h,k}(t)&:=&\xi^{t_0,y+k}(t)h-\xi^{t_0,y}(t)h-
\eta^{t_0,y}(t) (h,k)
\\[1pt]
&=&\int_{t_0}^te^{(t-s)A}DB
\bigl(s,Y^{t_0,y+k}(s) \bigr)\xi^{t_0,y+k}(s)h\udt s
\\[1pt]
&&{} -\int_{t_0}^te^{(t-s)A}DB
\bigl(s,Y^{t_0,y}(s) \bigr)\xi^{t_0,y}(s)h\udt s
\\[1pt]
&&{}-\int_{t_0}^te^{(t-s)A}D^2B
\bigl(s,Y^{t_0,y}(s) \bigr) \bigl(\xi^{t_0,y}(s)h,
\xi^{t_0,y}(s)k \bigr)\udt s
\\[1pt]
&&{} -\int_{t_0}^te^{(t-s)A}DB
\bigl(s,Y^{t_0,y}(s) \bigr)\eta^{t_0,y}(s) (h,k)\udt s
\\[1pt]
&=&\int_{t_0}^te^{(t-s)A}DB
\bigl(s,Y^{t_0,y+k}(s) \bigr)\xi^{t_0,y+k}(s)h\udt s
\\[1pt]
&&{} -\int_{t_0}^te^{(t-s)A}DB
\bigl(s,Y^{t_0,y}(s) \bigr)\xi^{t_0,y+k}(s)h\udt s
\\[1pt]
&&{}+\int_{t_0}^te^{(t-s)A}DB
\bigl(s,Y^{t_0,y}(s) \bigr)\xi^{t_0,y+k}(s)h\udt s
\\[1pt]
&&{} -\int_{t_0}^te^{(t-s)A}DB
\bigl(s,Y^{t_0,y}(s) \bigr)\xi^{t_0,y}(s)h\udt s
\\[1pt]
&&{} -\int_{t_0}^te^{(t-s)A}D^2B
\bigl(s,Y^{t_0,y}(s) \bigr) \bigl(\xi^{t_0,y}(s)h,
\xi^{t_0,y}(s)k \bigr)\udt s
\\[1pt]
&&{} -\int_{t_0}^te^{(t-s)A}DB
\bigl(s,Y^{t_0,y}(s) \bigr)\eta^{t_0,y}(s) (h,k)\udt s
\\[1pt]
&=&\int_{t_0}^te^{(t-s)A}DB
\bigl(s,Y^{t_0,y}(s) \bigr)\tilde r^{t_0,y,h,k}(s)\udt s
\\[1pt]
&&{} +\int_{t_0}^te^{(t-s)A}
\bigl[DB \bigl(s,Y^{t_0,y+k}(s) \bigr)-DB \bigl(s,Y^{t_0,y}(s) \bigr)
\bigr]\xi^{t_0,y+k}(s)h\udt s
\\[1pt]
&&{}-\int_{t_0}^te^{(t-s)A}D^2B
\bigl(s,Y^{t_0,y}(s) \bigr) \bigl(\xi^{t_0,y}(s)h,
\xi^{t_0,y}(s)k \bigr)\udt s
\\[1pt]
&=&\int_{t_0}^te^{(t-s)A}DB
\bigl(s,Y^{t_0,y}(s) \bigr)\tilde r^{t_0,y,h,k}(s)\udt s
\\[1pt]
&&{} +\int_{t_0}^te^{(t-s)A}
\biggl[\int_0^1D^2B \bigl(s,\alpha
Y^{t_0,y+k}(s)+(1-\alpha)Y^{t_0,y}(s) \bigr)\udt\alpha
\\[1pt]
&&{}\times \bigl(\xi^{t_0,y+k}(s)h,Y^{t_0,y+k}(s)-Y^{t_0,y}(s)
\bigr)
\\[1pt]
&&{}-D^2B \bigl(s,Y^{t_0,y}(s) \bigr) \bigl(\xi
^{t_0,y}(s)h,\xi^{t_0,y}(s)k \bigr) \biggr]\udt s
\\[1pt]
&=&\int_{t_0}^te^{(t-s)A}DB
\bigl(s,Y^{t_0,y}(s) \bigr)\tilde r^{t_0,y,h,k}(s)\udt s
\\[1pt]
&&{} +\int_{t_0}^te^{(t-s)A}
\biggl[\int_0^1D^2B \bigl(s,\alpha
Y^{t_0,y+k}(s)+(1-\alpha)Y^{t_0,y}(s) \bigr)\udt\alpha
\\[1pt]
&&{} -D^2B \bigl(s,Y^{t_0,y}(s) \bigr) \biggr]
\bigl(\xi^{t_0,y+k}(s)h,Y^{t_0,y+k}(s)-Y^{t_0,y}(s) \bigr)\udt
s
\\[1pt]
&&{}+\int_{t_0}^te^{(t-s)A}
\bigl[D^2B \bigl(s,Y^{t_0,s}(s) \bigr) \bigl(
\xi^{t_0,y+k}(s)h,Y^{t_0,y+k}(s)-Y^{t_0,y}(s) \bigr)
\\
&&{} -D^2B \bigl(s,Y^{t_0,y}(s) \bigr) \bigl(\xi
^{t_0,y}(s)h,\xi^{t_0,y}(s)k \bigr) \bigr]\udt s
\\
&=&\int_{t_0}^te^{(t-s)A}DB
\bigl(s,Y^{t_0,y}(s) \bigr)\tilde r^{t_0,y,h,k}(s)\udt s
\\
&&{} +\int_{t_0}^te^{(t-s)A}
\biggl[\int_0^1D^2B \bigl(s,\alpha
Y^{t_0,y+k}(s)+(1-\alpha)Y^{t_0,y}(s) \bigr)\udt\alpha
\\
&&{} -D^2B \bigl(s,Y^{t_0,y}(s) \bigr) \biggr]
\bigl(\xi^{t_0,y+k}(s)h,Y^{t_0,y+k}(s)-Y^{t_0,y}(s) \bigr)\udt
s
\\
&&{}+\int_{t_0}^te^{(t-s)A}D^2B
\bigl(s,Y^{t_0,y}(s) \bigr)
\\
&&{}\times \bigl[ \bigl(\xi^{t_0,y+k}(s)h,Y^{t_0,y+k}(s)-Y^{t_0,y}(s)
\bigr)- \bigl(\xi^{t_0,y+k}(s)h,\xi^{t_0,y}(s)k \bigr)
\\
&&{} + \bigl(\xi^{t_0,y+k}(s)h,\xi^{t_0,y}(s)k \bigr)-
\bigl(\xi^{t_0,y}(s)h,\xi^{t_0,y}(s)k \bigr) \bigr]\udt s
\\
&=&\int_{t_0}^te^{(t-s)A}DB
\bigl(s,Y^{t_0,y}(s) \bigr)\tilde r^{t_0,y,h,k}(s)\udt s
\\
&&{} +\int_{t_0}^te^{(t-s)A}
\biggl[\int_0^1D^2B \bigl(s,\alpha
Y^{t_0,y+k}(s)+(1-\alpha)Y^{t_0,y}(s) \bigr)\udt\alpha
\\
&&{} -D^2B \bigl(s,Y^{t_0,y}(s) \bigr) \biggr]
\bigl(\xi^{t_0,y+k}(s)h,Y^{t_0,y+k}(s)-Y^{t_0,y}(s) \bigr)\udt
s
\\
&&{}+\int_{t_0}^te^{(t-s)A}D^2B
\bigl(s,Y^{t_0,y}(s) \bigr)
\\
&&{}\times \bigl(\xi^{t_0,y+k}(s)h,Y^{t_0,y+k}(s)-Y^{t_0,y}(s)-
\xi^{t_0,y}(s)k \bigr)\udt s
\\
&&{} +\int_{t_0}^te^{(t-s)A}D^2B
\bigl(s,Y^{t_0,y}(s) \bigr)
\\
&&{}\times  \bigl(\xi^{t_0,y+k}(s)h-
\xi^{t_0,y}(s)h,\xi^{t_0,y}(s)k \bigr)\udt s.
\end{eqnarray*}
These calculations together with (\ref{eqstimaY}) and (\ref
{eqstimaDxi}) imply that
\begin{eqnarray*}
&& \bigl\llVert\tilde r^{t_0,y,h,k}(t)\bigr\rrVert_E
\\[-1pt]
&&\qquad \leq C{
\llVert DB\rrVert}_{\infty
}\int_{t_0}^t {
\bigl\llVert\tilde r^{t_0,y,h,k}(s)\bigr\rrVert}_E\udt s
\\[-1pt]
&&\quad\qquad{} +C\int_{t_0}^t\biggl\llVert\int
_0^1 D^2B \bigl(s,\alpha
Y^{t_0,y+k}(s)+(1-\alpha)Y^{t_0,y} \bigr)\udt\alpha
\\[-1pt]
&&\quad\qquad{}-D^2B
\bigl(s,Y^{t_0,y}(s) \bigr)\biggr\rrVert_{L(E,E;E)}{\bigl\llVert\xi
^{t_0,y+k}(s)h\bigr\rrVert}_E
\\[-1pt]
&&\quad\qquad{}\times {\bigl\llVert
Y^{t_0,y+k}(s)-Y^{t_0,y}(s)\bigr\rrVert}_E\udt s
\\[-1pt]
&&\quad\qquad{} +C{\bigl\llVert D^2B\bigr\rrVert}_\infty\int
_{t_0}^t{\bigl\llVert\xi^{t_0,y+k}(s)h
\bigr\rrVert}_E
\\[-1pt]
&&\quad\qquad{}\times {\bigl\llVert Y^{t_0,y+k}(s)-Y^{t_0,y}(s)-
\xi^{t_0,y}(s)k\bigr\rrVert}_E\udt s
\\[-1pt]
&&\quad\qquad{} +C{\bigl\llVert D^2B\bigr\rrVert}_\infty\int
_{t_0}^t{\bigl\llVert\xi^{t_0,y+k}(s)h-\xi
^{t_0,y}(s)h\bigr\rrVert}_E\cdot{\bigl\llVert
\xi^{t_0,y}(s)k\bigr\rrVert}_E\udt s
\\[-1pt]
&&\qquad \leq C{\llVert DB\rrVert}_{\infty}\int_{t_0}^t{
\bigl\llVert\tilde r^{t_0,y,h,k}(s)\bigr\rrVert}_E\udt s
+C\cdot C_\xi\widetilde C_Y \llVert h\rrVert
_E\llVert k\rrVert_E
\\[-1pt]
&&\quad\qquad{}\times \int
_{t_0}^t\biggl\llVert\int_0^1
D^2B \bigl(s,\alpha Y^{t_0,y+k}(s)+(1-\alpha)Y^{t_0,y}
\bigr)\udt\alpha
\\[-1pt]
&&\quad\qquad{}-D^2B \bigl(s,Y^{t_0,y}(s) \bigr)\biggr\rrVert
_{L(E,E;E)}\udt s
\\[-1pt]
&&\quad\qquad{}+C\cdot C_\xi\big\|D^2B\big\|_\infty\llVert h\rrVert_E
\\[-1pt]
&&\quad\qquad{}\times \int
_{t_0}^t\biggl\llVert\int_0^1
\xi^{t_0,\alpha(y+k)+(1-\alpha)y}(s)k\udt\alpha-\xi^{t_0,y}(s)k\biggr
\rrVert
_E\udt s
\\[-1pt]
&&\quad\qquad{}+C\cdot C_\xi\widetilde C_\xi T {\bigl
\llVert D^2B\bigr\rrVert}_\infty{\llVert h\rrVert
}_E{\llVert k\rrVert}_E^2
\\[-1pt]
&&\qquad \leq C{\llVert DB\rrVert}_{\infty}\int_{t_0}^t{
\bigl\llVert\tilde r^{t_0,y,h,k}(s)\bigr\rrVert}_E\udt
s+C_1{\llVert h\rrVert}_E{\llVert k\rrVert
}_E
\\[-1pt]
&&\quad\qquad{}\times \int_{t_0}^t\biggl\llVert
\int_0^1 D^2B \bigl(s,\alpha
Y^{t_0,y+k}(s)+(1-\alpha)Y^{t_0,y} \bigr)\udt\alpha
\\[-1pt]
&&\quad\qquad{}-D^2B \bigl(s,Y^{t_0,y}(s) \bigr)\biggr
\rrVert_{L(E,E;E)}\udt s
\\[-1pt]
&&\quad\qquad{}+C_2{\llVert h\rrVert}_E\int
_{t_0}^t\biggl\llVert\int_0^1
\xi^{t_0,y+\alpha k}(s)\udt\alpha-\xi^{t_0,y}(s)\biggr\rrVert
_{L(E,E)}\udt s{\llVert k\rrVert}_E
\\[-1pt]
&&\quad\qquad{} +C_3{\llVert h\rrVert}_E{\llVert k
\rrVert}_E^2.
\end{eqnarray*}
Finally, by an application of Gronwall's lemma
\begin{eqnarray*}
\frac{{\llVert \tilde r^{t_0,y,h,k}(t)\rrVert }_E}{{\llVert
k\rrVert }_E}&\leq& C_4{\llVert h\rrVert}_E
\\
&&{}\times \biggl[\int_{t_0}^t\biggl
\llVert\int_0^1 D^2B \bigl(s,
\alpha Y^{t_0,y+k}(s)+(1-\alpha)Y^{t_0,y} \bigr)\udt\alpha
\\
&&{}-D^2B \bigl(s,Y^{t_0,y}(s) \bigr)\biggr
\rrVert_{L(E,E;E)}\udt s
\\
&&{} +\int_{t_0}^t\biggl\llVert\int
_0^1\xi^{t_0,y+\alpha
k}(s)\udt\alpha-
\xi^{t_0,y}(s)\biggr\rrVert_{L(E,E)}\udt s + {\llVert k\rrVert
}_E \biggr]
\end{eqnarray*}
and such quantity goes to $0$ uniformly in ${\llVert h\rrVert }_E\leq N$
$\forall N>0$ when ${\llVert k\rrVert }_E$ goes to $0$ by Lebesgue's dominated
convergence theorem.

Our next step is to study the continuity of the second derivative
computed above. We have
\begin{eqnarray*}
&&\eta^{t_0,y}(t) (h,k)-\eta^{t_0,w}(t) (h,k)
\\
&&\qquad =\int_{t_0}^te^{(t-s)A}
\bigl[D^2B \bigl(s,Y^{t_0,y}(s) \bigr) \bigl(
\xi^{t_0,y}(s)h,\xi^{t_0,y}(s)k \bigr)
\\
&&\quad\qquad{} -D^2B \bigl(s,Y^{t_0,w}(s) \bigr) \bigl(\xi
^{t_0,w}(s)h,\xi^{t_0,w}(s)k \bigr) \bigr]\udt s
\\
&&\quad\qquad{}  +\int_{t_0}^te^{(t-s)A}
\bigl[DB \bigl(s,Y^{t_0,y}(s) \bigr)\eta^{t_0,y}(s) (h,k)
\\
&&\quad\qquad{} -DB \bigl(s,Y^{t_0,w}(s) \bigr)\eta^{t_0,w}(s)
(h,k) \bigr]\udt s
\\
&&\qquad = I_1+I_2;
\end{eqnarray*}
then
\begin{eqnarray*}
I_1&=&\int_{t_0}^te^{(t-s)A}
\bigl[D^2B \bigl(s,Y^{t_0,y}(s) \bigr) \bigl(\xi
^{t_0,y}(s)h,\xi^{t_0,y}(s)k \bigr)
\\
&&{}-D^2B \bigl(s,Y^{t_0,w}(s) \bigr) \bigl(\xi
^{t_0,y}(s)h,\xi^{t_0,y}(s)k \bigr)
\\
&&{}+D^2B \bigl(s,Y^{t_0,w}(s) \bigr) \bigl(\xi
^{t_0,y}(s)h,\xi^{t_0,y}k \bigr)
\\
&&{}-D^2B \bigl(s,Y^{t_0,w}(s) \bigr) \bigl(\xi
^{t_0,w}(s)h,\xi^{t_0,w}(s)k \bigr) \bigr]\udt s
\\
&=&\int_{t_0}^te^{(t-s)A}
\bigl[D^2B \bigl(s,Y^{t_0,y}(s) \bigr)-D^2B
\bigl(s,Y^{t_0,w}(s) \bigr) \bigr]
\bigl(\xi^{t_0,y}(s)h,\xi^{t_0,y}(s)k \bigr)\udt s
\\
&&{} +\int_{t_0}^te^{(t-s)A}D^2B
\bigl(s,Y^{t_0,w}(s) \bigr) \bigl( \bigl[\xi^{t_0,y}(s)-
\xi^{t_0,w} \bigr]h,\xi^{t_0,y}k \bigr)\udt s
\\
&&{}+\int_{t_0}^te^{(t-s)A}D^2B
\bigl(s,Y^{t_0,w}(s) \bigr) \bigl(\xi^{t_0,w}(s)h, \bigl[
\xi^{t_0,y}(s)-\xi^{t_0,w}(s) \bigr]k \bigr)\udt s
\end{eqnarray*}
and
\begin{eqnarray*}
I_2&=&\int_{t_0}^te^{(t-s)A}DB
\bigl(s,Y^{t_0,y}(s) \bigr) \bigl[\eta^{t_0,y}(s) (h,k)-
\eta^{t_0,w}(s) (h,k) \bigr]\udt s
\\
&&{}+\int_{t_0}^te^{(t-s)A} \bigl[DB
\bigl(s,Y^{t_0,y}(s) \bigr)-DB \bigl(s,Y^{t_0,w}(s) \bigr) \bigr]
\eta^{t_0,w}(s) (h,k)\udt s.
\end{eqnarray*}
Recalling all the estimates previously obtained and the fact that both
${\llVert Y^{t_0,y}(t)\rrVert }_E$ and ${\llVert \xi
^{t_0,y}(t)\rrVert
}_{L(E,E)}$ are
bounded uniformly in $t$, denoting with $C_H$ the H\"older constant of
$D^2B$, we get
\begin{eqnarray*}
&& \bigl\llVert \eta^{t_0,y}(t) (h,k)-\eta^{t_0,w}(t) (h,k)\bigr
\rrVert_E
\\
&&\qquad \leq C\cdot C_H\int_{t_0}^t\bigl
\llVert Y^{t_0,y}(s)-Y^{t_0,w}(s)\bigr\rrVert_E^\alpha
\bigl\llVert\xi^{t_0,y}(s)h\bigr\rrVert_E\bigl\llVert\xi
^{t_0,y}(s)k\bigr\rrVert_E\udt s
\\
&&\quad\qquad{}  +C{\bigl\llVert D^2B\bigr\rrVert}_\infty
\int_{t_0}^t\bigl\llVert\xi^{t_0,y}(s)-
\xi^{t_0,w}(s)\bigr\rrVert_{L(E,E)}^\alpha\bigl\llVert
\xi^{t_0,y}(s)-\xi^{t_0,w}(s)\bigr\rrVert_{L(E,E)}^{1-\alpha}
\\
&&\quad\qquad{}\times \bigl[{\llVert h\rrVert}_E{\bigl\llVert\xi
^{t_0,y}(s)k\bigr\rrVert}_E+\bigl\llVert
\xi^{t_0,w}(s)h\bigr\rrVert_E{\llVert k\rrVert
}_E \bigr]\udt s
\\
&&\quad\qquad{}+C\llVert DB\rrVert_\infty\int_{t_0}^t
\bigl\llVert\eta^{t_0,y}(s) (h,k)-\eta^{t_0,w}(s) (h,k)\bigr
\rrVert_E\udt s
\\
&&\quad\qquad{} +C\llVert DB\rrVert_\infty\int_{t_0}^t
\bigl\llVert Y^{t_0,y}(s)-Y^{t_0,w}(s)\bigr\rrVert
_E^\alpha\bigl\llVert Y^{t_0,y}(s)-Y^{t_0,w}(s)
\bigr\rrVert_E^{1-\alpha}
\\
&&\quad\qquad{} \times \bigl\llVert\eta^{t_0,w}(s) (h,k)\bigr\rrVert_E
\udt s
\\
&&\qquad \leq C_5{\llVert h\rrVert}_E{\llVert k\rrVert
}_E\llVert y-w\rrVert_E^\alpha
+C_6\int_{t_0}^t\bigl\llVert
\eta^{t_0,y}(s) (h,k)-\eta^{t_0,w}(s) (h,k)\bigr\rrVert
_E\udt s
\end{eqnarray*}
hence
\[
\bigl\llVert\eta^{t_0,y}(t) (h,k)-\eta^{t_0,w}(t) (h,k)\bigr
\rrVert_E\leq C_7{\llVert h\rrVert}_E{
\llVert k\rrVert}_E{\llVert y-w\rrVert}_E^\alpha
\]
which shows that the second Fr\'echet derivative of the map $y\mapsto
Y^{t_0,y}(t)$ is $\alpha$-H\"older continuous.\vspace*{6pt}

\textit{Continuity with respect to the initial time}.
Fix $t\in[0,T]$, $\omega\in\Omega_0$ (that we do not write, as
before) and $\epsilon>0$ and consider two initial times $s_1$ and
$s_2$, with $s_1<s_2$ for simplicity. Since we assume that $y\in\sL
^p$ or $y\in\aC$, we can find $\delta$ such that
\begin{eqnarray*}
&& \bigl\llVert Y^{s_2,y}(t)-Y^{s_1,y}(t)\bigr\rrVert
_E
\\
&&\qquad \leq\biggl\llVert e^{(t-s_2)A} \bigl(1-e^{(s_2-s_1)A} \bigr)y
\\
&&\quad\qquad{} +\int_{s_2}^te^{(t-r)A}
\bigl[B \bigl(r,Y^{s_2,y}(r) \bigr)-B \bigl(r,Y^{s_1,y}(r) \bigr)
\bigr]\udt r
\\
&&\quad\qquad{} -\int_{s_1}^{s_2}e^{(t-r)A}B
\bigl(r,Y^{s_1,y}(r) \bigr)\udt r-\int_{s_1}^{s_2}e^{(t-r)A}
\Sigma\udt W(r)\biggr\rrVert_E
\\
&&\qquad\leq C\bigl\llVert\bigl(1-e^{ (s_2-s_1 )A} \bigr)y\bigr\rrVert
_E+C\llVert DB\rrVert_\infty\int_{s_2}^t
\bigl\llVert Y^{s_2,y}(r)-Y^{s_1,y}(r)\bigr\rrVert
_E\udt r
\\
&&\quad\qquad{} +C\llVert B\rrVert_\infty\llvert s_2-s_1
\rrvert+C\llVert\Sigma\rrVert_\infty\bigl\llvert W(s_2)-W(s_1)
\bigr\rrvert
\\
&&\qquad\leq C\llVert DB\rrVert_\infty\int_{s_2}^t
\bigl\llVert Y^{s_2,y}(r)-Y^{s_1,y}(r)\bigr\rrVert
_E\udt r+C\epsilon
\end{eqnarray*}
for $\llvert s_2-s_1\rrvert <\delta$, because $e^{sA}$ is strongly
continuous and $W(\cdot,\omega)$ is continuous.

The conclusion follows using Gronwall's lemma, $\epsilon$ being arbitrary.
\end{appendix}


\section*{Acknowledgments}
We are grateful to the referees for their careful comments and
stimulating questions that lead us to improve the exposition and revise
several details for the final version of our work.


%

\printaddresses
\end{document}